\documentclass[12pt]{article}
\usepackage{amsmath,amsfonts,amssymb}
\usepackage[enableskew]{youngtab}
\usepackage{theorem}
\usepackage[usenames]{color}
\nonstopmode

\newtheorem{theorem}{Theorem}[section]
\newtheorem{proposition}[theorem]{Proposition}
\newtheorem{lemma}[theorem]{Lemma}
\newtheorem{definition}[theorem]{Definition}
\newtheorem{corollary}[theorem]{Corollary}
\newtheorem{conjecture}[theorem]{Conjecture}
\newtheorem{claim}[theorem]{Claim}
\theorembodyfont{\rmfamily}
\newtheorem{example}[theorem]{Example}
\newtheorem{remark}[theorem]{Remark}

\begin{document}
$\,$\vspace{10mm}

\begin{center}
{\textsf {\huge Paths and Kostka--Macdonald Polynomials}}
\vspace{10mm}\\
{\textsf{\LARGE  ${}^{\mbox{\small a}}$Anatol N. Kirillov and
${}^{\mbox{\small b}}$Reiho Sakamoto}}
\vspace{20mm}\\
{\textsf {${}^{\mbox{\small{a}}}$Research Institute for Mathematical Sciences,}}
\vspace{-1mm}\\
{\textsf {Kyoto University, Sakyo-ku,}}
\vspace{-1mm}\\
{\textsf {Kyoto, 606-8502, Japan}}
\vspace{-1mm}\\
{\textsf {kirillov@kurims.kyoto-u.ac.jp}}
\vspace{3mm}\\
{\textsf {${}^{\mbox{\small{b}}}$Department of Physics, Graduate School of Science,}}
\vspace{-1mm}\\
{\textsf {University of Tokyo, Hongo, Bunkyo-ku,}}
\vspace{-1mm}\\
{\textsf {Tokyo, 113-0033, Japan}}
\vspace{-1mm}\\
{\textsf {reiho@rs.kagu.tus.ac.jp}}
\vspace{40mm}
\end{center}

\begin{abstract}
\noindent
We give several equivalent combinatorial
descriptions of the space of states for the 
box-ball systems, and connect certain partition
functions for these models with the 
$q$-weight multiplicities of the tensor product of
the fundamental representations of the Lie algebra $\mathfrak{gl}(n)$.
As an application, we give an elementary proof of the 
special case $t=1$ of the Haglund--Haiman--Loehr formula.
Also, we propose a new class of combinatorial
statistics that naturally generalize the
so-called energy statistics.\bigskip\\
Mathematics Subject Classification (2000) 05E10,~20C35.\\
Key words and phrases:~Crystals, Paths, Energy and Tau functions, Box--Ball 
systems, Kostka--Macdonald polynomials. 
\end{abstract}

\pagebreak

\section{Introduction}
The purpose of the present paper is two-fold.
First of all, we would like to give an introduction to
the beautiful Combinatorics related with Box-Ball Systems, and secondly,
to relate the latter with the ``Classical Combinatorics''
revolving around transportation matrices, tabloids, the 
Lascoux--Sch\"{u}tzenberger  statistics charge, Macdonald polynomials,
\cite{M},\cite{Sag},
Haglund--Haiman--Loehr's formula \cite{HHL}, and so on. 
As a result of our investigations, we will prove that two statistics
naturally appearing in the context of Box-Ball systems,
namely {\it energy function} and {\it  tau-function},
have nice combinatorial properties.
More precisely, the statistics energy $E$ 
is an example of a generalized machonian 
statistics \cite{Kir}, Section~2, whereas the statistics $\tau$ related
with Kostka--Macdonald polynomials, see Section~5.2 of the present paper.

Box-Ball Systems (BBS for short) were invented by
Takahashi--Satsuma \cite{TS,Tak} as a wide class of discrete
integrable soliton systems.
In the simplest case, BBS are described by simple combinatorial
procedures using box and balls.
Despite its simple outlook, it is known that
the BBS have various remarkably deep properties;
\begin{itemize}
\item
Time evolution of the BBS coincides with isomorphism
of the crystal bases \cite{HHIKTT,FOY}.
Thus the BBS possesses quantum integrability.
\item
BBS are ultradiscrete (or tropical) limit of the
usual soliton systems \cite{TTMS,KSY}.
Thus the BBS possesses classical integrability
at the same time.
\item
Inverse scattering formalism of the BBS coincides with
the rigged configuration bijection originating in
completeness problem of the Bethe states \cite{KOSTY,Sak3}.
\end{itemize}
Let us say a few words about the main results of our paper.
\begin{itemize}
\item
In the case of statistics {\it tau},
our main result can be formulated as a computation
of the corresponding  partition function for the BBS in terms of the
values of the Kostka--Macdonald polynomials at $t=1.$
\item
In the case of the statistics {\it energy},
our result can be formulated as 
an interpretation of the corresponding
partition function for the BBS as the $q$-weight 
multiplicity in the tensor product of
the fundamental representations of the Lie 
algebra $\mathfrak{gl}(n).$ We {\it expect}
that the same statement is valid for the BBS 
corresponding to the tensor product of rectangular representations.

~~~~~~We are reminded that a {\it $q$-analogue of the multiplicity} of a highest 
weight $\lambda$ in the tensor product $\bigotimes_{a=1}^{L}~V_{s_{a} \omega_{r_a}}$
of the highest weight $s_a~\omega_{r_a},$~$a=1,\ldots,L,$ irreducible representations 
$V_{s_{a} \omega_{r_{a}}}$ of the Lie algebra $\mathfrak {gl}(n)$  is defined as
$$q\mbox{-Mult }[V_{\lambda}~ :~ \bigotimes_{a=1}^{L}~V_{s_{a} \omega_{r_{a}}}] =
\sum_{\eta}~K_{\eta, R}~K_{\eta, \lambda}(q),$$
where  $K_{\eta, R}$ stands for the parabolic Kostka number corresponding to the 
sequence of rectangles $R:= \{ (s_{a}^{r_{a}}) \}_{a=1,\ldots,L},$ 
see e.g. \cite{Kir}, \cite{KiSh}.
\item
We give several equivalent descriptions of
{\it paths} which appear in the 
description of partition functions for BBS:
in terms of transportation matrices, 
tabloids, plane partitions.
We expect that such interpretations may be
helpful for better understanding connections of
the BBS and other integrable models such as melting crystals \cite{ORV}, 
$q$-difference Toda lattices \cite{GLO}, ... .
\end{itemize} 

Our result about connections of the energy partition
functions for BBS and $q$-weight 
multiplicities suggests a deep hidden connections
between partition functions for the BBS 
and characters of the Demazure modules, solutions to the
$q$-difference Toda equations, cf.\cite{GLO}. 

As an interesting open problem we want to give raise
a question about an interpretation 
of the sums $\sum_{\eta}~K_{\eta,R}~K_{\eta,\lambda}(q,t),$ where
 $K_{\eta,\lambda}(q,t)$ 
denotes the Kostka--Macdonald polynomials \cite{M},
as {\it refined partition functions} for the BBS 
corresponding to the tensor product of rectangular representations 
$R=\{(s_a^{r_a}) \}_{1 \le a \le n}$.
See Conjecture \ref{conj:tau}. In other words, one can ask: what is a meaning of the 
second statistics (see \cite{HHL}) in the Kashiwara theory \cite{Kas}  of crystal 
bases (of type A) ?

Organization of the present paper is as follows.
In Section 2, we review necessary facts from the
Kirillov--Reshetikhin crystals.
Especially we explain an explicit algorithm
to compute the combinatorial $R$
and the energy function.
In section 3, we introduce several combinatorial
descriptions of paths.
Then we define several statistics on paths
such as Haglund's statistics, energy statistics $\bar{E}$
and tau statistics $\tau^{r,s}$.
In Section 4, we collect necessary facts from
the BBS which will be used in the next section.
In Section 5, we present our main result
(Theorem \ref{th:main}) as well as several
relating conjectures.
We conjecture that $\tau^{r,s}$ gives independent
statistics depending on one parameter $r$
although they all give rise to the
unique generating function up to constant shift of power.
In Section 6, we show that the energy statistics
$\bar{E}$ belong to the class of statistics $\tau^{r,s}$
(Theorem \ref{th:E=tau}).
Therefore $\tau^{r,s}$ gives a natural extension
of the energy statistics $\bar{E}$.

\section{Kirillov--Reshetikhin crystal}
\subsection{$A^{(1)}_n$ type crystal}
Let $W_s^{(r)}$ be a $U_q(\mathfrak{g})$ Kirillov--Reshetikhin
module, where we shall consider the case $\mathfrak{g}=A_n^{(1)}$.
The module $W_s^{(r)}$ is indexed by a Dynkin node
$r\in I=\{1,2,\ldots,n\}$ and $s\in\mathbb{Z}_{>0}$.
As a $U_q(A_n)$-module, $W_s^{(r)}$ is isomorphic to the irreducible
module corresponding to the partition $(s^r)$.
For arbitrary $r$ and $s$, the module $W_s^{(r)}$
is known to have crystal bases \cite{Kas,KMN2}, 
which we denote by $B^{r,s}$.
As the set, $B^{r,s}$ is consisting of all column strict
semi-standard Young tableaux
of depth $r$ and width $s$ over the alphabet $\{1,2,\ldots,n+1\}$.

For the algebra $A_n$, let $P$ be the weight lattice,
$\{\Lambda_i\in P|i\in I\}$ be the fundamental roots,
$\{\alpha_i\in P|i\in I\}$ be the simple roots,
and $\{h_i\in \mathrm{Hom}_\mathbb{Z}(P,\mathbb{Z})|i\in I\}$
be the simple coroots.
As a type $A_n$ crystal, $B=B^{r,s}$ is equipped with
the Kashiwara operators $e_i,f_i:B\longrightarrow B\cup\{0\}$
and $\mathrm{wt}:B\longrightarrow P$ ($i\in I$) satisfying
\begin{eqnarray*}
&&f_i(b)=b'\Longleftrightarrow e_i(b')=b\quad \mbox{ if }b,b'\in B,\\
&&\mathrm{wt}\left(f_i(b)\right)=\mathrm{wt}(b)-\alpha_i
\quad \mbox{ if }f_i(b)\in B,\\
&&\langle h_i,\mathrm{wt}(b)\rangle =\varphi_i(b)-\varepsilon_i(b).
\end{eqnarray*}
Here $\langle\cdot ,\cdot\rangle$ is the natural pairing
and we set 
$\varepsilon_i(b)=\max\{m\ge0\mid \tilde{e}_i^m b\ne0\}$
and $\varphi_i(b)=\max\{m\ge0\mid \tilde{f}_i^m b\ne0\}$.
Actions of the Kashiwara operators
$\tilde{e}_i$, $\tilde{f}_i$
for $i\in I$ coincide with the one described in \cite{KN}.
Since we do not use explicit forms of these operators,
we omit the details.
See \cite{O} for complements of this section.
Note that in our case $A_n$, we have $P=\mathbb{Z}^{n+1}$
and $\alpha_i=\epsilon_i-\epsilon_{i+1}$ where
$\epsilon_i$ is the $i$-th canonical unit vector of $\mathbb{Z}^{n+1}$.
We also remark that $\mathrm{wt}(b)=(\lambda_1,\cdots,\lambda_{n+1})$
is the weight of $b$, i.e., $\lambda_i$ counts the number of letters
$i$ contained in tableau $b$.

For two crystals $B$ and $B'$,
one can define the tensor product
$B\otimes B'=\{b\otimes b'\mid b\in B,b'\in B'\}$.
The actions of the Kashiwara operators on tensor
product have simple form.
Namely, the operators 
$\tilde{e}_i,\tilde{f}_i$ act on $B\otimes B'$ by
\begin{eqnarray*}
\tilde{e}_i(b\otimes b')&=&\left\{
\begin{array}{ll}
\tilde{e}_i b\otimes b'&\mbox{ if }\varphi_i(b)\ge\varepsilon_i(b')\\
b\otimes \tilde{e}_i b'&\mbox{ if }\varphi_i(b) < \varepsilon_i(b'),
\end{array}\right. \\
\tilde{f}_i(b\otimes b')&=&\left\{
\begin{array}{ll}
\tilde{f}_i b\otimes b'&\mbox{ if }\varphi_i(b) > \varepsilon_i(b')\\
b\otimes \tilde{f}_i b'&\mbox{ if }\varphi_i(b)\le\varepsilon_i(b'),
\end{array}\right.
\end{eqnarray*}
and $\mathrm{wt}(b\otimes b')=\mathrm{wt}(b)+
\mathrm{wt}(b')$.
We assume that $0\otimes b'$ and $b\otimes 0$ as $0$.
Then it is known that there is a unique crystal isomorphism
$R:B^{r,s}\otimes B^{r',s'}
\stackrel{\sim}{\rightarrow}B^{r',s'}\otimes B^{r,s}$.
We call this map (classical) combinatorial $R$
and usually write the map $R$ simply by $\simeq$.

Let us consider the affinization of the crystal $B$.
As the set, it is
\begin{equation}
\mathrm{Aff}(B)=\{b[d]\, |\, b\in B,\, d\in\mathbb{Z}\}.
\end{equation}
There is also explicit algorithm for
actions of the affine Kashiwara operators $\tilde{e}_0$, $\tilde{f}_0$
in terms of  the promotion operator \cite{Shimo}.
For the tensor product
$b[d]\otimes b'[d']\in
\mathrm{Aff}(B)\otimes\mathrm{Aff}(B')$,
we can lift the (classical)
combinatorial $R$ to affine case
as follows:
\begin{equation}\label{eq:affineR}
b[d]\otimes b'[d']\stackrel{R}{\simeq}
\tilde{b}'[d'-H(b\otimes b')]\otimes
\tilde{b}[d+H(b\otimes b')],
\end{equation}
where $b\otimes b'\simeq \tilde{b}'\otimes \tilde{b}$
is the isomorphism of (classical) combinatorial $R$.
The function $H(b\otimes b')$ is called the energy function.
We will give explicit forms of the combinatorial $R$ and
energy function in the next section.

\subsection{Combinatorial $R$ and energy function}
We give explicit description of the combinatorial $R$-matrix
(combinatorial $R$ for short)
and energy function on $B^{r,s}\otimes B^{r',s'}$.
To begin with we define few terminologies about Young tableaux.
Denote rows of a Young tableaux $Y$ by $y_1,y_2,\ldots y_r$
from the top to bottom.
Then row word $row(Y)$ is defined by concatenating rows as
$row(Y)=y_ry_{r-1}\ldots y_1$.
Let $x=(x_1,x_2,\ldots )$ and $y=(y_1,y_2,\ldots )$ be two partitions.
We define concatenation of $x$ and $y$ by the partition
$(x_1+y_1,x_2+y_2,\ldots )$.

\begin{proposition}[\cite{Shimo}]\label{prop:shimozono}
$b\otimes b'\in B^{r,s}\otimes B^{r',s'}$ is mapped to
$\tilde{b}'\otimes \tilde{b}\in B^{r',s'}\otimes B^{r,s}$
under the combinatorial $R$, i.e.,
\begin{equation}
b\otimes b'\stackrel{R}{\simeq}\tilde{b}'\otimes\tilde{b},
\end{equation}
if and only if
\begin{equation}
(b'\leftarrow row(b))=(\tilde{b}\leftarrow row(\tilde{b}')).
\end{equation}
Moreover, the energy function $H(b\otimes b')$ is given by
the number of nodes of $(b'\leftarrow row(b))$
outside the concatenation of partitions
$(s^r)$ and $({{s'}^{r'}})$.
\hfill$\square$
\end{proposition}

We define another normalization of the energy function
$\bar{H}$ such that for $b\otimes b'\in B^{r,s}\otimes
B^{r',s'}$,
\begin{equation}
\bar{H}(b\otimes b):=\min (r,r')\cdot\min (s,s')-
H(b\otimes b).
\end{equation}
For special cases of $B^{1,s}\otimes B^{1,s'}$,
the function $H$ is called unwinding number and
$\bar{H}$ is called winding number in \cite{NY}.
Explicit values for the case
$b\otimes b'\in B^{1,1}\otimes B^{1,1}$
are given by
\begin{equation}
H(b\otimes b')=
\chi (b<b'),\qquad
\bar{H}(b\otimes b')=
\chi (b\geq b'),
\end{equation}
where $\chi(\mathrm{True})=1$ and
$\chi(\mathrm{False})=0$.

In order to describe the
algorithm for finding $\tilde{b}$ and $\tilde{b}'$ from
the data $(b'\leftarrow row(b))$,
we introduce a terminology.
Let $Y$ be a tableau, and $Y'$ be a subset of $Y$
such that $Y'$ is also a tableau.
Consider the set theoretic subtraction $\theta =Y\setminus Y'$.
If the number of nodes contained in $\theta$ is $r$
and if the number of nodes of $\theta$ contained in
each row is always 0 or 1,
then $\theta$ is called vertical $r$-strip.

Given a tableau
$Y=(b'\leftarrow row(b))$, let $Y'$ be the upper left
part of $Y$ whose shape is $(s^r)$.
We assign numbers from 1 to $r's'$
for each node contained in $\theta =Y\setminus Y'$
by the following procedure.
Let $\theta_1$ be the vertical $r'$-strip of $\theta$
as upper as possible.
For each node in $\theta_1$,
we assign numbers 1 through $r'$ from the bottom to top.
Next we consider $\theta\setminus\theta_1$,
and find the vertical $r'$ strip $\theta_2$
by the same way.
Continue this procedure until all nodes of $\theta$
are assigned numbers up to $r's'$.
Then we apply inverse bumping procedure according
to the labeling of nodes in $\theta$.
Denote by $u_1$ the integer which is ejected
when we apply inverse bumping procedure starting
from the node with label 1.
Denote by $Y_1$ the tableau such that
$(Y_1\leftarrow u_1)=Y$.
Next we apply inverse bumping procedure
starting from the node of $Y_1$ labeled by 2,
and obtain the integer $u_2$ and tableau $Y_2$.
We do this procedure until we obtain $u_{r's'}$
and $Y_{r's'}$.
Finally, we have
\begin{equation}
\tilde{b}'=(\emptyset\leftarrow u_{r's'}u_{r's'-1}\cdots
u_1),\qquad
\tilde{b}=Y_{r's'}.
\end{equation}\begin{example}
Consider the following tensor product:
$$b\otimes b'=
\Yvcentermath1
\young(114,235)\otimes\young(23,34,45)\in
B^{2,3}\otimes B^{3,2}.$$
{}From $b$, we have $row(b)=235114$, hence we have
$$\left
(
\Yvcentermath1\young(23,34,45)
\leftarrow 235114
\right)=
\newcommand{\yonsan}{4_3}
\newcommand{\sanni}{3_2}
\newcommand{\yonichi}{4_1}
\newcommand{\goyon}{5_4}
\newcommand{\yongo}{4_5}
\newcommand{\sanroku}{3_6}
\Yvcentermath1
\young(113\yonsan ,225,\sanroku\sanni ,\yongo\yonichi ,\goyon) .
$$
Here subscripts of each node indicate the order of
inverse bumping procedure.
For example, we start from the node $4_1$ and obtain
$$\left
(
\Yvcentermath1
\young(1234,235,34,4,5)
\leftarrow 1
\right)=
\Yvcentermath1
\young(1134,225,33,44,5),
\qquad\mathrm{therefore,}\qquad
Y_1=
\newcommand{\yonsan}{4_3}
\newcommand{\yonni}{4_2}
\newcommand{\goyon}{5_4}
\newcommand{\yongo}{4_5}
\newcommand{\sanroku}{3_6}
\Yvcentermath1
\young(123\yonsan ,235,\sanroku\yonni ,\yongo ,\goyon)
,\qquad
u_1=1.
$$
Next we start from the node $4_2$ of $Y_1$.
Continuing in this way, we obtain
$u_6u_5\cdots u_1=321421$ and
$Y_6=\Yvcentermath1
\young(334,455)$.
Since $(\emptyset\leftarrow 321421)=
\Yvcentermath1\young(11,22,34)$,
we obtain
$$\Yvcentermath1
\young(114,235)\otimes\young(23,34,45)\simeq
\young(11,22,34)\otimes\young(334,455)\, ,
\qquad
H\left
(
\Yvcentermath1
\young(114,235)\otimes\young(23,34,45)
\right)
=3.
$$
Note that the energy function is derived from the
concatenation of shapes of $b$ and $b'$,
i.e., $\Yvcentermath1\yng(5,5,2)\,$.
\hfill$\square$
\end{example}

\section{Combinatorics on the set of paths}
\subsection{Combinatorics}
\subsubsection{Transportation matrices and tabloids}
Let $n$ be a positive integer,
$\alpha=(\alpha_1,\cdots,\alpha_n)$ and
$\beta=(\beta_1,\cdots.\beta_n)$ be two compositions
of the same size.
Denote by $M_n(\alpha,\beta)$ the set of matrices
$M=(m_{i,j})_{1\leq i,j\leq n}$ such that
\begin{equation}
m_{i,j}\in\mathbb{Z}_{\geq 0},\quad
\sum_j m_{i,j}=\alpha_i,\quad
\sum_i m_{i,j}=\beta_j.
\end{equation}
Remind that a tabloid of shape $\alpha$ and weight $\beta$
is a filling of the shape $\alpha$ by the
numbers $1,2,\cdots,n$ in such a way that
the number of $i$'s appearing in the filling
is equal to $\beta_i$.
It is clear that the number of tabloids
of shape $\alpha$ and weight $\beta$ is equal to the multinomial coefficient
$$\frac{(\beta_1+\beta_2+\cdots +\beta_n)!}
{\beta_1!\beta_2!\cdots\beta_n!}.$$

A row (column) weakly strict tabloid of shape
$\alpha$ and weight $\beta$ is a filling of the shape
$\alpha$ by numbers $1,2,\cdots,n$ such that the numbers
along each row (column) are weakly increasing and
$\beta_i$ is equal to the number of $i$'s appearing
in the filling.

\begin{example}
Take
$\alpha =(2,1,3,1)$,
$\beta =(1,3,0,2,1)$, then
$$\Yvcentermath1
\young(45,1,224,2)\,\, .$$
is a row weakly strict tabloid of shape $\alpha$ and weight $\beta$.
\hfill$\square$
\end{example}

We denote by $\mathrm{Tab}(\alpha,\beta)$ the set of all
row weakly strict tabloids of shape $\alpha$ and weight $\beta$.
It is well-known that there exists a bijection between
the sets $M_n(\alpha,\beta)$ and $\mathrm{Tab}(\alpha,\beta)$.
Namely, given a matrix $m=(m_{ij})\in M_n(\alpha ,\beta)$,
we fill the row $\alpha_i$ of the shape $\alpha$ by the numbers
$1^{m_{i1}}, 2^{m_{i2}},\cdots, n^{m_{in}}$.
For example, let
$$\left(
\begin{array}{lllll}
0& 0& 0& 1& 1\\
1& 0& 0& 0& 0\\
0& 2& 0& 1& 0\\
0& 1& 0& 0& 0
\end{array}
\right)\,\, .$$
Then the corresponding row weakly strict tabloid is
$\Yvcentermath1
\young(45,1,224,2)\,\, .$
To each tabloid $T$, one can associate the reading word,
namely the word obtained by reading the filling of tabloid $T$
from the right to the left starting from the top row.
For example, for the tabloid $T$ displayed above,
$w(T)=5414222$.

If weight $\mu$ of a tabloid $T$ is a partition,
we define the charge $c(T)$ of tabloid $T$ to be the charge
$c(w(T))$ of the reading word.
See page 242 of \cite{M} for the definition of the
Lascoux--Sch\"{u}tzenberger charge \cite{LS}. 

\begin{example}
Take standard tabloid
$$T=\Yvcentermath1
\young(45,1,367,2)\,\, ,$$
then
$$w(T)=5_3 4_2 1_0 7_4 6_3 3_1 2_0$$
and therefore
$c(T)=3+2+0+4+3+1+0=13$.
\hfill$\square$
\end{example}

\subsubsection{Plane partitions}

Let $\lambda$ be a partition.
{\it A plane partition of shape $\lambda$} is a 
tabloid $\pi$ of shape $\lambda$ such that
the numbers in each row and each column 
are weakly decreasing. For example,
$$ \pi= \Yvcentermath1
\young(754,744,33,3)\,\, ,$$
is a plane partitions of shape $(3,3,2,1).$

A plane partition $\pi$ has a three-dimensional diagram,
consisting of the points 
$(i,j,k)$ with integer coordinates
such that $(i,j) \in \lambda$ and  
$1 \le k \le \pi(i,j),$ where $\pi(i,j) \in \pi$
is the number that is located in 
the box $(i,j) \in \lambda.$
By definition, the size $|\pi|$ of a plane partition 
$\pi$ of shape $\lambda$ is
$|\pi|=\sum_{(i,j) \in \lambda}~\pi(i,j).$

Let $\alpha$ and $\beta$ be two compositions of the same size. Denote by 
$\mathcal{PP}(\alpha,\beta)$ the set of plane partitions $\pi$ such that
$$\sum_{k \ge 0}~\pi(i,i+k)=\sum_{j \ge k} \alpha_j,~~~
\sum_{k \ge 0}~\pi(i+k,i)=\sum_{j \ge k} \beta_j.$$ 
Finally, let us remind two classical results
\begin{enumerate}
\item[(A)]
(P. MacMahon, see e.g. \cite{M}, page 81)
Let $l,m,n$ be three positive integers, and $B$ be the box with 
side-lengths $l.m,n.$ Then
$$ \sum_{\pi \subset B}~q^{|\pi|} =
\prod_{(i,j,k) \in B}~{1-q^{i+j+k-1} \over 1-q^{i+j+k-2}}.$$
\item[(B)]
(Robinson--Schensted--Knuth, see e.g.
\cite{Sag}, Chapter~3) There are bijections
$$\mathcal{M}(\alpha,\beta)
\leftarrow\hspace{-3.5mm}
\xrightarrow{\,\, 1:1\,\,}
\mathcal{PP}(\alpha, \beta)
\leftarrow\hspace{-3.5mm}
\xrightarrow{\,\, 1:1\,\,}
\mathrm{Tab}(\alpha,\beta).$$
\end{enumerate}

\subsection{Paths}
\subsubsection{$B^{1,1}$ type paths}
Let $\alpha$ be a partition of size $n$.
A path $p$ of type $B^{1,1}$ and weight
$\alpha=(\alpha_1,\alpha_2,\cdots,\alpha_n)$
is a sequence of positive integers $a_1a_2\cdots a_n$
such that $\alpha_i=\#\{ j|a_j=i\}$.
We denote by $\mathcal{P}(\alpha)$ the set of
all paths of type $B^{1,1}$ and weight $\alpha$.
A path $p$ is called a highest weight path if
the sequence $a_1a_2\cdots a_n$ satisfies
the Yamanouchi condition.
We denote by $\mathcal{P}_+(\alpha)$ the set of
all $B^{1,1}$ type highest paths with
weight $\alpha$.
It is well known that the total number of $B^{1,1}$
type paths of weight $\alpha=(\alpha_1,\alpha_2,
\cdots,\alpha_n)$ is equal to the multinomial
coefficient
$$\frac{(\alpha_1+\alpha_2+\cdots +\alpha_n)!}
{\alpha_1!\alpha_2!\cdots\alpha_n!},$$
and there are  bijections
$$\mathcal{P}(\alpha)
\leftarrow\hspace{-3.5mm}
\xrightarrow{\,\, 1:1\,\,}
\mathrm{Mat}_n(\alpha,1^n)
\leftarrow\hspace{-3.5mm}
\xrightarrow{\,\, 1:1\,\,}
\mathrm{Tab}(\alpha,1^n).$$
Let us describe the general
prescription to get the corresponding tabloid
from a given path.
Let the path $a_1a_2\cdots a_n\in\mathcal{P}(\alpha)$,
we recursively add letters to the tabloid
according to $a_1$, $a_2$, $\cdots$, $a_n$ as follows.
Starting from the empty tabloid,
assume that we have done up to $a_{i-1}$
and have gotten a tabloid $T^{(i-1)}$.
Then we add the letter $i$ to the right
of the $a_i$-th row of $T^{(i-1)}$
and get $T^{(i)}$.
For example, the path $p=4221343$ can be related to
the following transportation matrix and row strict tabloid
$$M=\left(
\begin{array}{lllllll}
0& 0& 0& 1& 0& 0& 0\\
0& 1& 1& 0& 0& 0& 0\\
0& 0& 0& 0& 1& 0& 1\\
1& 0& 0& 0& 0& 1& 0
\end{array}
\right)\,\, ,~~~~~~~
T=\Yvcentermath1
\young(4,23,57,16)\,\, .$$

\subsubsection{General rectangular paths}
More generally, we define path to be an arbitrary
element of tensor product of crystals $B^{r_1,s_1}\otimes
B^{r_2,s_2}\otimes\cdots\otimes B^{r_L,s_L}$.
Recall that for type $A^{(1)}_n$ case,
$B^{r,s}$ is, as the set, consisting of semi-standard
tableaux over alphabet $\{1,2,\cdots n+1\}$,
and tensor product of crystals $B\otimes B'$ is, as a set,
cartesian product of two sets $B$ and $B'$.
Crystal graph structure on the set $B\otimes B'$
is is given according to \cite{KN}.
Weight $\lambda=(\lambda_1,\cdots,\lambda_{n+1})$
of a path $b=b_1\otimes b_2\otimes\cdots\otimes b_L\in
B^{r_1,s_1}\otimes
B^{r_2,s_2}\otimes\cdots\otimes B^{r_L,s_L}$ is
given by 
\begin{equation}
\lambda_i=\mbox{total number of letters }i
\mbox{ contained in tableaux }B^{r_1,s_1},\cdots, B^{r_L,s_L}.
\end{equation}
For example,
$$\Yvcentermath1
\young(12,24,35)\otimes\young(12,33)\otimes \young(11,45)\in
B^{3,2}\otimes B^{2,2}\otimes B^{2,2}$$
is a path of rectangular shape
$R=((2^3),(2^2),(2^2))$,
and its weight is $\lambda =(4,3,3,2,2)$.
Note that the number of standard (i.e., weight of $(1^N)$)
rectangular shape $R=\{(r_a^{s_a})_{a=1,2,\cdots ,L}\}$ paths
is equal to the generalized multinomial coefficient
\begin{equation}
\genfrac{[}{]}{0pt}{}{N}{R_1,\cdots,R_L}
=\frac{N!}{\prod_a H_{R_a}(1)},
\end{equation}
where $N=\sum_a r_as_a$, and
for any diagram $\lambda$, $H_\lambda (q)$
denotes the hook polynomial (see definition on page 45 of \cite{M})
corresponding to diagram $\lambda$.

\paragraph{Comments.}
Summarizing, one has the following (equivalent)
combinatorial descriptions of the set of 
(crystal) paths of type
$B^{1,s_1}\otimes B^{1,s_2}\otimes\cdots\otimes B^{1,s_n}$ and 
weight $\alpha=(\alpha_1, \ldots,\alpha_L)$ as the set of 

$(a)$ transportation matrices $\mathcal{M}_{L}(\alpha ,s),$

$(b)$ row weakly increasing tabloids $\mathcal{T}(\alpha ,s),$

$(c)$ plane partitions $\mathcal{PP}(\alpha ,s).$

\noindent
For the given path
$b_1\otimes b_2\otimes\cdots\otimes b_n\in
B^{1,s_1}\otimes B^{1,s_2}\otimes\cdots\otimes B^{1,s_n}$,
the corresponding element in
$\mathcal{T}(\alpha ,s)$ is determined
as follows.
Staring from the empty tabloid,
we recursively add letters to the tabloid
according to $b_1$, $b_2$, $\cdots$, $b_n$
as follows.
Assume we have done up to $b_{i-1}$
and have gotten the tabloid $T^{(i-1)}$.
Denote the number of $k$ contained in
$b_i$ by $x_k$.
Then, for all $k$,
we add letters $i$ for $x_k$ times
to the right of the $k$-th row of $T^{(i-1)}$
and get $T^{(i)}$.

\begin{example}
Consider the path
$$\Yvcentermath1
\young(33)\otimes\young(1123)\otimes\young(13)\otimes\young(233)$$
of type $B^{1,2}\otimes B^{1,4}\otimes B^{1,2}\otimes B^{1,3}$
and weight $\alpha =(3,2,6)$.
The corresponding tabloid and transportation matrix are
$$T=\Yvcentermath1\young(223,24,112344)\,\, ,\qquad
M=\left(
\begin{array}{llll}
0& 2& 1& 0\\
0& 1& 0& 1\\
2& 1& 1& 2
\end{array}
\right)\, .$$
To find the plane partition which corresponds to
the tabloid $T$ (or matrix $M$), one can apply
the Robinson--Schensted--Knuth algorithm \cite{Knu}
to the multi-permutation
$$w:=\left(
\begin{array}{c}
1\, 1\, 1\, 2\, 2\, 3\, 3\, 3\, 3\, 3\, 3\\
2\, 2\, 3\, 2\, 4\, 1\, 1\, 2\, 3\, 4\, 4
\end{array}
\right)$$
which corresponds to the tabloid $T$. One has
$$w\leftarrow\hspace{-3.5mm}
\xrightarrow{\,\, RSK\,\,}
\left(\,
\Yvcentermath1\young(1122344,224,3)\,\, ,\,\,
\young(1112333,233,3)\,
\right)\, .$$
Finally, the plane partition we are looking
for, can be obtained from the pair of semi-standard
Young tableaux displayed above by gluing the
Gelfand--Tsetlin patterns that correspond
to the Young tableaux in question:
\begin{center}
\unitlength 13pt
\begin{picture}(25,7.5)(0,-2.5)
\put(0,2){$
\left(
\begin{array}{c}
7\quad 3\quad 1\quad 0\\
5\quad 2\quad 1\\
4\quad 2\\
2
\end{array}\,\,\, ,\,\,\,
\begin{array}{c}
7\quad 3\quad 1\quad 0\\
7\quad 3\quad 1\\
4\quad 1\\
3
\end{array}
\right)\,\,
\longleftrightarrow
$}
\put(16,2){
$\begin{array}{cccc}
7& 7& 4& 3\\
5& 3& 3& 1\\
4& 2& 1& 1\\
2& 2& 1& 0
\end{array}$
}
\multiput(17,0.4)(0.2,-0.2){10}{\circle*{0.1}}
\multiput(17,1.6)(0.2,-0.2){16}{\circle*{0.1}}
\multiput(17,2.8)(0.2,-0.2){22}{\circle*{0.1}}
\multiput(17,4.0)(0.2,-0.2){28}{\circle*{0.1}}
\multiput(18.2,4.0)(0.2,-0.2){22}{\circle*{0.1}}
\multiput(19.4,4.0)(0.2,-0.2){16}{\circle*{0.1}}
\multiput(20.6,4.0)(0.2,-0.2){10}{\circle*{0.1}}
\put(19,-2.2){2}
\put(20.2,-2.2){6}
\put(21.4,-2.2){8}
\put(22.7,-2.2){11}
\put(22.7,-1.0){11}
\put(22.7,0.2){5}
\put(22.7,1.4){3}
\end{picture}
\end{center}
The result is a plane partition from the
set $\mathcal{PP}((2,4,2,3),(3,2,6,0))$.
\hfill$\square$
\end{example}
\begin{remark}
For the reader's convenience,
let us recall the way to get the corresponding
Gelfand--Tsetlin pattern from a given
semi-standard tableaux.
By looking contents of a semi-standard
tableau $T$,
one can define a sequence of partitions
\begin{equation}\label{def:GT}
\emptyset =\mu^{(0)}\subset
\mu^{(1)}\subset\cdots\subset\mu^{(n)} =\mu
\end{equation}
such that each skew diagram
$\mu^{(i)}\setminus\mu^{(i-1)}$
$(1\leq i\leq n)$ is a horizontal
strip, see e.g., Chapter I of \cite{M}.
Starting from the sequence of partitions
(\ref{def:GT}), one can define the
corresponding Gelfand--Tsetlin pattern
$x:=x(T)$ by the following rule
\begin{equation}
x^{(i)}(T)=\mbox{shape}(\mu^{(i)}),
\qquad (1\leq i\leq n).
\end{equation}
It is known that thus obtained $x$
indeed satisfies the defining
properties of the
Gelfand--Tsetlin patterns.
\hfill$\square$
\end{remark}
\begin{remark}
One of the basic properties of the
BBS is that the second Young tableau\footnote{
Equivalently the upper part of the
corresponding plane partition.}
(of weight $\beta$) obtained by means
of the Robinson--Schensted--Knuth algorithm,
is conserved under the dynamics of the BBS
\cite{Fuk} (see also \cite{TTS,Ari} for the other
connections between the simplest BBS and the RSK algorithm).
Nowadays, conserved quantities and linearization
parameters (or angle variables) of the BBS
are completely
determined in the most general settings \cite{KOSTY},
and, surprisingly enough, they are elegantly
described by the so-called  (unrestricted)
rigged configurations  \cite{KKR,KR,KI,KSS,Sch,DS}.
The latter result is a consequence of a deep theorem
stated in Lemma 8.5 of \cite{KSS}.
\hfill$\square$
\end{remark}

\subsection{Statistics on the set of paths}

\subsubsection{Energy statistics}
For a path $b_1\otimes b_2\otimes
\cdots\otimes b_L\in B^{r_1,s_1}\otimes
B^{r_2,s_2}\otimes\cdots\otimes B^{r_L,s_L}$,
let us define elements
$b_j^{(i)}\in B^{r_j,s_j}$ for $i<j$
by the following isomorphisms of the combinatorial $R$;
\begin{eqnarray}
&&b_1\otimes b_2\otimes\cdots\otimes
b_{i-1}\otimes b_{i}\otimes\cdots\otimes
b_{j-1}\otimes b_j\otimes\cdots\nonumber\\
&\simeq&
b_1\otimes b_2\otimes\cdots\otimes
b_{i-1}\otimes b_{i}\otimes\cdots\otimes
b_j^{(j-1)}\otimes b_{j-1}'\otimes\cdots
\nonumber\\
&\simeq&\cdots\nonumber\\
&\simeq&
b_1\otimes b_2\otimes\cdots\otimes
b_{i-1}\otimes b_{j}^{(i)}\otimes
\cdots\otimes
b_{j-2}'\otimes b_{j-1}'\otimes\cdots,
\end{eqnarray}
where we have written
$b_k\otimes b_j^{(k+1)}\simeq
b_j^{(k)}\otimes b_k'$
assuming that $b_j^{(j)}=b_j$.

For a given path  $p=b_1\otimes b_2\otimes
\cdots\otimes b_L\in B^{r_1,s_1}\otimes
B^{r_2,s_2}\otimes\cdots\otimes B^{r_L,s_L}$,
define statistics $\bar{E}(p)$ by
\begin{equation}
\bar{E}(p)=\sum_{i<j}
\bar{H}(b_i\otimes b_j^{(i+1)}).
\end{equation}

Define the statistics $\mathrm{maj}(p)$ by
\begin{equation}\label{def:maj}
\mathrm{maj}(p)=
\sum_{i<j}H(b_i\otimes b_j^{(i+1)}).
\end{equation}
For example, consider a path
$a=a_1\otimes a_2\otimes\cdots\otimes a_L
\in (B^{1,1})^{\otimes L}$.
In this case, we have $a_j^{(i)}=a_i$,
since the combinatorial $R$ act on
$B^{1,1}\otimes B^{1,1}$ as identity.
Therefore, we have
\begin{equation}
\mathrm{maj}(a)
=\sum_{i=1}^{L-1}(L-i)\chi (a_i<a_{i+1}).
\end{equation}

Define another statistics {\it tau} as follows.
\begin{definition}\label{def:generaltau}
For the path $p\in B^{r_1,s_1}\otimes
B^{r_2,s_2}\otimes\cdots\otimes B^{r_L,s_L}$,
define $\tau^{r,s}$ by
\begin{equation}
\tau^{r,s}(p)=\mathrm{maj}
(u^{(r)}_s\otimes p),
\end{equation}
where $u^{(r)}_s$ is the highest element of $B^{r,s}$.
\hfill$\square$
\end{definition}
We use abbreviation $\tau$ for the statistics
$\tau^{1,1}$ on $B^{1,1}$ type paths
$a\in (B^{1,1})^{\otimes L}$, i.e.,
\begin{equation}\label{eq:taumu}
\tau (a)=\mathrm{maj}(1\otimes a)=\mathrm{maj}(a)+L~(1-\delta_{1,a_{1}}),
\end{equation}
where $a_1$ denotes the first letter of the path $a$.
This $\tau$ is a special case of the tau functions
for the box-ball systems \cite{KSY,Sak1}
which originate as ultradiscrete limit of
the tau functions for the KP hierarchy \cite{JM}.
\begin{definition}\label{def:taumu}
For composition $\mu=(\mu_1,\mu_2,\cdots,\mu_n)$,
write $\mu_{[i]}=\sum_{j=1}^i \mu_j$
with convention $\mu_{[0]}=0$.
Then we define a generalization
of $\tau$ by
\begin{equation}
\tau_\mu(a)=\sum_{i=1}^{n}
\tau(a_{[i]}),
\end{equation}
where
\begin{equation}
a_{[i]}=
a_{\mu_{[i-1]}+1}\otimes
a_{\mu_{[i-1]}+2}\otimes\cdots\otimes
a_{\mu_{[i]}}
\in (B^{1,1})^{\otimes\mu_i}.
\end{equation}
\hfill$\square$
\end{definition}
Note that we have
$a=a_{[1]}\otimes a_{[2]}\otimes\cdots\otimes
a_{[n]}$, i.e., the path $a$ is partitioned
according to $\mu$.
It is convenient to identify $\tau_\mu$
as statistics on a tabloid of shape $\mu$ whose
reading word coincides with
the partitioned path according to $\mu$.
For example, to path
$p={\it abcdefgh}$
and the composition $\mu=(3,2,3)$
one associates the tabloid
$$\Yvcentermath1
\young(cba,ed,hgf)\,\, .$$

\subsubsection{Statistics charge}
Let $p$ be a path of type $B^{1,1},$ denote by $T(p)$ the corresponding row strict 
tabloid and by $w(T(p))$ its reading word. Define the charge of a path $p$ to be the 
charge of tabloid $T(p),$ i.e. the charge of the reading word $w(T(p)).$ For example, 
take  $p=4221343.$  Then $w(T(p))=4327561,$ and therefore, $c(p)=3+2+1+4+3+3+0=16.$

If $\mu$ is a composition, define $c_{\mu}(p)=\sum_{i}~c(p_{[i]}),$ where 
$$p_{[i]}= p_{\mu_{[i-1]}+1}~p_{\mu_{[i-1]}+2}~\cdots~p_{\mu_{[i]}},$$
cf. Definition~3.3.
\begin{lemma} One has
$$ \tau_{\mu}(p)+c_{\mu}(p)=\sum_{i}~\biggl( {m_{i}+1 \choose 2}-
\mu_{i}~\delta_{1,p_{1}^{[i]}} \biggr),$$
where $p_1^{[i]}$ denotes the first letter of the path $p_{[i]}.$
\hfill$\square$
\end{lemma}
Proof follows from two simple observations that
$$\tau(p)+c(1 \otimes p)={L \choose 2},
\qquad
c(1 \otimes p) -c(p)= L~\delta_{1,p_{1}},$$
where $L$ denotes the length of path $p.$

\paragraph{Comments.}
It follows from \cite{KSY} that on the set of semi-standard
Young tableaux, i.e., on the set of highest weight paths,
the statistics tau coincides with statistics {\it cocharge}.
Therefore, one can consider the statistics tau as a natural
extension of the statistics cocharge from the set of
semi-standard tableaux to the set of tabloids,
or on the set of transportation matrices.
\hfill$\square$

\subsubsection{Haglund's statistics}
\paragraph{Tableaux language description}
For a given path
$a=a_1\otimes a_2\otimes\cdots\otimes a_L
\in (B^{1,1})^{\otimes L}$, associate tabloid $t$
of shape $\mu$
whose reading word coincides with $a$.
This correspondence is the same as those
used in Definition \ref{def:taumu}.
Denote the cell at the $i$-th row,
$j$-th column
(we denote the coordinate by $(i,j)$) of
the tabloid $t$ by $t_{ij}$.
Attacking region of the cell at $(i,j)$
is all cells
$(i,k)$ with $k<j$ or $(i+1,k)$ with $k>j$.
In the following diagram, gray zonal
regions are the attacking regions of
the cell $(i,j)$.

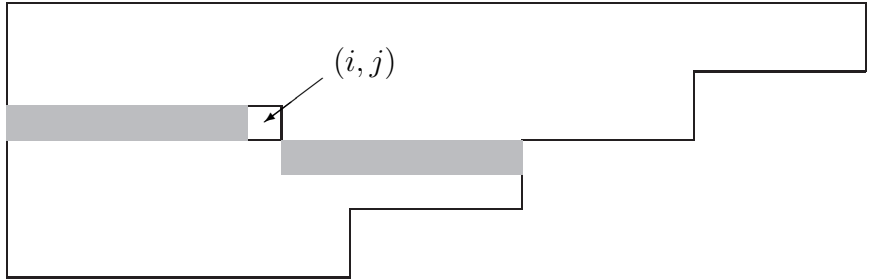
\unitlength 13pt
\begin{picture}(10,9)(-2,1.5)
\put(0,2){\line(0,1){8}}
\put(0,2){\line(1,0){5}}
\put(0,10){\line(1,0){25}}
\multiput(5,2)(5,2){4}{\line(1,0){5}}
\multiput(10,2)(5,2){4}{\line(0,1){2}}
\multiput(7,6)(0,1){2}{\line(1,0){1}}
\multiput(7,6)(1,0){2}{\line(0,1){1}}
\color[cmyk]{0,0,0,0.3}
\put(0,6){\rule{91pt}{13pt}}
\put(8,5){\rule{91pt}{13pt}}
\color{black}
\put(9.2,7.8){\vector(-4,-3){1.7}}
\put(9.5,8){$(i,j)$}
\end{picture}

\noindent
Follow \cite{HHL}, define $|\mathrm{Inv}_{ij}|$ by
\begin{equation}
|\mathrm{Inv}_{ij}|=
\#\{(k,l)\in \mbox{ attacking region for } (i,j)\, |\,
t_{kl}>t_{ij}\}.
\end{equation}
Then we define
\begin{equation}
|\mathrm{Inv}_\mu (a)|=\sum_{(i,j)\in \mu}
|\mathrm{Inv}_{ij}|.
\end{equation}

If we have $t_{(i-1)j}<t_{ij}$,
then the cell $(i,j)$ is called by {\it descent}.
Then define
\begin{equation}
{\rm Des}_\mu(a)=
\sum_{{\rm all\, descent\, }(i,j)}
(\mu_i-j).
\end{equation}
Note that $(\mu_i-j)$ is the arm length
of the cell $(i,j)$.

\paragraph{Path language description}
Consider two paths
$a^{(1)},a^{(2)}\in (B^{1,1})^{\otimes\mu}$.
We denote by
$a^{(1)}\otimes a^{(2)}=a_1\otimes a_2
\otimes\cdots\otimes a_{2\mu}$.
Then we define
\begin{equation}
\mathrm{Inv}_{(\mu ,\mu)}(a^{(1)},a^{(2)})=
\sum_{k=1}^{\mu}\sum_{i=k+1}^{k+\mu -1}
\chi (a_k<a_i).
\end{equation}
For more general cases
$a^{(1)}\in (B^{1,1})^{\otimes\mu_1}$ and
$a^{(2)}\in (B^{1,1})^{\otimes\mu_2}$
satisfying $\mu_1>\mu_2$, we define
\begin{equation}
\mathrm{Inv}_{(\mu_1,\mu_2)}(a^{(1)},a^{(2)}):=
\mathrm{Inv}_{(\mu_1,\mu_1)}
{(a^{(1)},1^{\otimes (\mu_1-\mu_2)}
\otimes a^{(2)})}.
\end{equation}
Then the above definition of $|\mathrm{Inv}_\mu(a)|$
is equivalent to
\begin{equation}
|\mathrm{Inv}_\mu(a)|=
\sum_{i=1}^{n-1}\mathrm{Inv}_{(\mu_i,\mu_{i+1})}.
\end{equation}

For example, consider the following tabloid
($a=2312133212$);
$$t=\Yvcentermath1\young(312132,2123),\qquad
|\mathrm{Inv}_{(6,4)}(a)|=10.$$
We associate the paths
$a^{(1)}=231213$ and $a^{(2)}=1^{\otimes 2}3212$.
Then
\begin{align*}
|\mathrm{Inv}_{(a^{(1)},a^{(2)})}|=&
\chi (a_1<a_2)+\chi (a_1<a_3)+\chi (a_1<a_4)+\chi (a_1<a_5)+\chi (a_1<a_6)+\\
&
\chi (a_2<a_3)+\chi (a_2<a_4)+\chi (a_2<a_5)+\chi (a_2<a_6)+\chi (a_2<a_7)+\\
&
\chi (a_3<a_4)+\chi (a_3<a_5)+\chi (a_3<a_6)+\chi (a_3<a_7)+\chi (a_3<a_8)+\\
&
\chi (a_4<a_5)+\chi (a_4<a_6)+\chi (a_4<a_7)+\chi (a_4<a_8)+\chi (a_4<a_9)+\\
&
\chi (a_5<a_6)+\chi (a_5<a_7)+\chi (a_5<a_8)+\chi (a_5<a_9)+\chi (a_5<a_{10})+\\
&
\chi (a_6<a_7)+\chi (a_6<a_8)+\chi (a_6<a_9)+\chi (a_6<a_{10})+\chi (a_6<a_{11})\\
=&(1+0+1+0+1)+(0+0+0+0+0)+
(1+0+1+0+0)+\\
&(0+1+0+0+1)+(1+0+0+1+1)+(0+0+0+0+0)=10.
\end{align*}

Consider two paths
$a^{(1)}\in (B^{1,1})^{\otimes\mu_1}$ and
$a^{(2)}\in (B^{1,1})^{\otimes\mu_2}$
satisfying $\mu_1\geq \mu_2$.
Denote $a=a^{(1)}\otimes a^{(2)}$.
Then define
\begin{equation}
\mathrm{Des}_{(\mu_1,\mu_2)}(a)=
\sum_{k=\mu_1-\mu_2+1}^{\mu_1}
(k-(\mu_1-\mu_2)-1)\chi (a_k<a_{k+\mu_2}).
\end{equation}
For the tableau $T$ of shape $\mu$ corresponding to
the path $a$, we define
\begin{equation}
\mathrm{Des}_\mu (T)=
\sum_{i=1}^n\mathrm{Des}_{(\mu_i,\mu_{i+1})}
(a_{[i]}\otimes a_{[i+1]}).
\end{equation}

\begin{definition}[\cite{Hag}]
For a path $a$,
statistics $\mathrm{maj}_\mu$ is defined by
\begin{equation}
\mathrm{maj}_\mu (a)=
\sum_{i=1}^{\mu_1}\mathrm{maj}(t_{1,i}\otimes t_{2,i}
\otimes\cdots\otimes t_{\mu'_i,i}).
\end{equation}
and $\mathrm{inv}_{\mu}(a)$ is defined by
\begin{equation}
\mathrm{inv}_\mu (a)=|\mathrm{Inv}_\mu(a)|-
\mathrm{Des}_\mu (a).
\end{equation}
\hfill$\square$
\end{definition}

If we associate to a given path $p \in {\cal P}(\lambda)$
with the shape $\mu$ tabloid $T$, we sometimes write
$\mathrm{maj}_\mu (p)=\mathrm{maj}(T)$ and
$\mathrm{inv}_\mu (p)=\mathrm{inv}(T)$.
\begin{example}\label{ex:222_42}
For highest weight paths of weight $\lambda =(2,2,2)$
and shape $\mu=(4,2)$,
the following is a list of the corresponding tabloids
associated with
data in the form $(\mathrm{maj}_\mu,\mathrm{inv}_\mu):$
$$\begin{array}{|c|c|c|}\hline
\Yvcentermath1\young(2211,33)\,\,\,\,(2,4)&
\Yvcentermath1\young(2121,33)\,\,\,\,(2,3)&
\Yvcentermath1\young(3121,32)\,\,\,\,(1,5)
\begin{picture}(0,1.7)\end{picture}
\\[4mm]\hline
\Yvcentermath1\young(1321,32)\,\,\,\,(1,4)&
\Yvcentermath1\young(3211,32)\,\,\,\,(0,6)&
\begin{picture}(0,1.7)\end{picture}
\\[4mm]\hline
\end{array}$$

Let us observe that
\begin{align*}
\sum_{p\in\mathcal{P}_+(\lambda)}
q^{\mathrm{inv}_{\mu} (p)}
t^{\mathrm{maj}_{\mu} (p)}=
q^6+q^4 t+ q^5 t+q^3 t^2+q^4 t^2
\end{align*}
which is different from
\begin{align*}
\tilde{K}_{\lambda ,\mu}(q,t)=
q^6+q^4 t+ q^5 t+q^2 t^2+q^4 t^2.
\end{align*}

Another interesting choice is
$\lambda =(2,2,2)$ and $\mu=(3,3)$.
The following is a list of all such paths
with corresponding statistics:
$$\begin{array}{|c|c|c|}\hline
\Yvcentermath1\young(211,332)\,\,\,\,(3,3)&
\Yvcentermath1\young(211,323)\,\,\,\,(3,2)&
\Yvcentermath1\young(121,332)\,\,\,\,(3,2)
\begin{picture}(0,1.7)\end{picture}
\\[4mm]\hline
\Yvcentermath1\young(121,323)\,\,\,\,(2,3)&
\Yvcentermath1\young(321,321)\,\,\,\,(0,6)&
\begin{picture}(0,1.7)\end{picture}
\\[4mm]\hline
\end{array}$$
Thus we have
\begin{align*}
\sum_{p\in\mathcal{P}_+(\lambda)}
q^{\mathrm{inv}_{\mu} (p)}
t^{\mathrm{maj}_{\mu} (p)}=
q^6+q^3 t^3+q^3 t^2+2 q^2 t^3
\end{align*}
which is again different from
\begin{align*}
\tilde{K}_{\lambda ,\mu}(q,t)=
q^6+q^4 t^2+q^3 t^3+q^3 t^2+q^2 t^2.
\end{align*}
\hfill$\square$
\end{example}

\section{Box-ball system}
In this section, we summarize basic facts about
the box-ball system which will be used in
the next section.
For our purpose, it is convenient to
express the isomorphism of the combinatorial
$R$
\begin{equation}
a\otimes b\simeq b'\otimes a'
\end{equation}
by the following vertex diagram:
\begin{center}
\unitlength 13pt
\begin{picture}(4,4)
\put(0,2.0){\line(1,0){3.2}}
\put(1.6,1.0){\line(0,1){2}}
\put(-0.6,1.8){$a$}
\put(1.4,0){$b'$}
\put(1.4,3.2){$b$}
\put(3.4,1.8){$a'$}
\put(4.1,1.7){.}
\end{picture}
\end{center}
Successive applications of the combinatorial $R$
is depicted by concatenating these vertices.

Following \cite{HHIKTT,FOY}, we define time evolution
of the box-ball system $T^{(a)}_l$.
Let $u_{l,0}^{(a)}=u_{l}^{(a)}\in B^{a,l}$ be the highest
element and $b_i\in B^{r_i,s_i}$.
Here the highest element
$u_{l}^{(a)}\in B^{a,l}$ is the tableau
whose $i$-th row is occupied by integers $i$.
For example, $u_4^{(3)}=\Yvcentermath1
\young(1111,2222,3333)\,$.
Define $u_{l,j}^{(a)}$ and $b_i'\in B^{r_i,s_i}$
by the following diagram.
\begin{equation}\label{def:T_l}
\unitlength 13pt
\begin{picture}(22,5)(0,-0.5)
\multiput(0,0)(5.8,0){2}{
\put(0,2.0){\line(1,0){4}}
\put(2,0){\line(0,1){4}}
}
\put(-1.4,1.8){$u_{l,0}^{(a)}$}
\put(1.7,4.2){$b_1$}
\put(1.7,-0.8){$b_1'$}
\put(4.2,1.8){$u_{l,1}^{(a)}$}
\put(7.5,4.2){$b_2$}
\put(7.5,-0.8){$b_2'$}
\put(10.0,1.8){$u_{l,2}^{(a)}$}
\multiput(11.5,1.8)(0.3,0){10}{$\cdot$}
\put(14.7,1.8){$u_{l,L-1}^{(a)}$}
\put(17,0){
\put(0,2.0){\line(1,0){4}}
\put(2,0){\line(0,1){4}}
}
\put(18.7,4.2){$b_L$}
\put(18.7,-0.8){$b_L'$}
\put(21.2,1.8){$u_{l,L}^{(a)}$}
\end{picture}
\end{equation}
$u_{l,j}^{(a)}$ are usually called {\it carrier}
and we set $u_{l,0}^{(a)}:=u_{l}^{(a)}$.
Then we define operator $T_l^{(a)}$ by
\begin{equation}
T_l^{(a)}(b)=b'=
b_1'\otimes b_2'\otimes\cdots\otimes b_L'.
\end{equation}
Recently \cite{Sak2,Sak3}, operators $T^{(a)}_l$ have used
to derive crystal theoretical meaning of the
rigged configuration bijection.

It is known (\cite{KOSTY} Theorem 2.7) that there
exists some $l\in\mathbb{Z}_{>0}$ such that
\begin{equation}
T^{(a)}_{l}=T^{(a)}_{l+1}=T^{(a)}_{l+2}=\cdots
(=:T^{(a)}_\infty ).
\end{equation}
If the corresponding path is $b\in (B^{1,1})^{\otimes L}$,
we have the following combinatorial description
of the box-ball system \cite{TS,Tak}.
We regard $\fbox{1}\in B^{1,1}$ as an empty box of
capacity 1, and $\fbox{$i$}\in B^{1,1}$
as a ball of label (or internal degree of freedom) $i$
contained in the box.
Then we have:
 
\begin{proposition}[\cite{HHIKTT}]\label{prop:bbs}
For a path $b\in (B^{1,1})^{\otimes L}$ of type $A^{(1)}_n$,
$T^{(1)}_\infty (b)$ is given by
the following procedure.
\begin{enumerate}
\item
Move every ball only once.

\item
Move the leftmost ball with label $n+1$
to the nearest right empty box.

\item
Move the leftmost ball with label $n+1$
among the rest to its nearest right
empty box.

\item
Repeat this procedure until all of the balls
with label $n+1$ are moved.

\item
Do the same procedure 2--4 for the balls with
label $n$.

\item
Repeat this procedure successively until all
of the balls with label $2$ are moved.
\end{enumerate}
\hfill$\square$
\end{proposition}
There are extensions \cite{HKT2,HKT3} of this box and ball algorithm
corresponding to generalizations of the box-ball
systems with respect to each affine Lie algebra \cite{HKT1,HKOTY}.
Using this box and ball interpretation,
our statistics $\tau (b)$ admits the following
interpretation.

\begin{theorem}[\cite{KSY} Theorem 7.4]\label{th:tau=rho}
For a path $b\in (B^{1,1})^{\otimes L}$ of type $A^{(1)}_n$,
$\tau (b)$ coincides with number of all balls
$2,\cdots,n+1$ contained in paths
$b$, $T^{(1)}_\infty (b)$, $\cdots$,
$( T^{(1)}_\infty )^{L-1} (b)$.
\hfill$\square$
\end{theorem}

\begin{example}
Consider the path $p=a\otimes b$ where
$a=4312111$, $b=4321111$.
We compute $\tau_{(7,7)}(p)$ in two ways.

\noindent (i)
First we compute by Eq.(\ref{eq:taumu}).
\begin{align*}
&\tau (a)=\mathrm{maj}(1\otimes a)=
7\cdot 1 +6\cdot 0+5\cdot 0+4\cdot 1+
3\cdot 0+2\cdot 0+1\cdot 0  =11,\\
&\tau (b)=\mathrm{maj}(1\otimes b)=
7\cdot 1 +6\cdot 0+5\cdot 0+4\cdot 0+
3\cdot 0+2\cdot 0+1\cdot 0  =7.
\end{align*}
Thus we obtain
$\tau_{(7,7)}(p)=\tau (a)+\tau (b)=11+7=18$.

\noindent (ii)
Next we use Theorem \ref{th:tau=rho}.
According to Proposition \ref{prop:bbs},
the time evolutions of the paths $a$ and
$b$ are as follows:
\begin{center}
\fbox{\hspace{-3pt}
$\begin{array}{lllllll}
 4 & 3 & 1 & 2 & 1 & 1 & 1 \\
 1 & 1 & 4 & 1 & 3 & 2 & 1 \\
 1 & 1 & 1 & 4 & 1 & 1 & 3 \\
 1 & 1 & 1 & 1 & 4 & 1 & 1 \\
 1 & 1 & 1 & 1 & 1 & 4 & 1 \\
 1 & 1 & 1 & 1 & 1 & 1 & 4 \\
 1 & 1 & 1 & 1 & 1 & 1 & 1
\end{array}$}\qquad\qquad
\fbox{\hspace{-3pt}
$\begin{array}{lllllll}
 4 & 3 & 2 & 1 & 1 & 1 & 1 \\
 1 & 1 & 1 & 4 & 3 & 2 & 1 \\
 1 & 1 & 1 & 1 & 1 & 1 & 4 \\
 1 & 1 & 1 & 1 & 1 & 1 & 1 \\
 1 & 1 & 1 & 1 & 1 & 1 & 1 \\
 1 & 1 & 1 & 1 & 1 & 1 & 1 \\
 1 & 1 & 1 & 1 & 1 & 1 & 1
\end{array}
$}
\end{center}
Here the left and right tables correspond to
$a$ and $b$, respectively.
Rows of left (resp. right) table represent
$a$, $T^{(1)}_\infty (a)$, $\cdots$,
$( T^{(1)}_\infty )^L (a)$ (resp., those for $b$)
from top to bottom.
Note that we omit all frames of tableaux of $B^{1,1}$
and symbols for tensor product.
Counting letters 2, 3 and 4 in each table,
we have $\tau (a)=11$, $\tau (b)=7$
and again we get $\tau_{(7,7)}(p)=11+7=18$.
Meanings of the above two dynamics
corresponding to paths $a$ and $b$ are
summarized as follows:
\begin{enumerate}
\item[$(a)$]
Dynamics of the path $a$.
In the first row, there are two solitons
(length two soliton $43$
and length one soliton 2),
and in the second row, there are also two solitons
(length one soliton 4 and length two soliton 32).
This is scattering of two solitons.
After the scattering, soliton 4 propagates
at velocity one and soliton 32 propagates
at velocity two without scattering.
\item[$(b)$]
Dynamics of the path $b$.
This shows free propagation of one soliton
of length three 432 at velocity three.
\end{enumerate}
\hfill$\square$
\end{example}

\section{Main results}
\subsection{Haglund--Haiman--Loehr formula}
Let $\tilde{H}_\mu (x;q,t)$ be the (integral form ) modified Macdonald
polynomials where $x$ stands for
infinitely many variables $x_1,x_2,\cdots$.
Here $\tilde{H}_\mu (x;q,t)$ is obtained by
simple plethystic substitution (see, e.g.,
section 2 of \cite{Hai}) from the original
definition of the Macdonald polynomials \cite{M}.
Schur function expansion of $\tilde{H}_\mu (x;q,t)$
is given by
\begin{equation}
\tilde{H}_\mu (x;q,t)=
\sum_\lambda\tilde{K}_{\lambda ,\mu}(q,t)
s_\lambda (x),
\end{equation}
where $\tilde{K}_{\lambda ,\mu}(q,t)$
stands for the following transformation of the
Kostka--Macdonald polynomials:
\begin{equation}
\tilde{K}_{\lambda ,\mu}(q,t)=
t^{n(\mu)}K_{\lambda ,\mu}(q,t^{-1}).
\end{equation}
Here we have used notation
$n(\mu)=\sum_i(i-1)\mu_i$.
Then the celebrated
Haglund--Haiman--Loehr (HHL) formula is
as follows.

\begin{theorem}[\cite{HHL}]
Let $\sigma:\mu\rightarrow\mathbb{Z}_{>0}$
be the filling of the Young diagram $\mu$
by positive integers $\mathbb{Z}_{>0}$,
and define $x^\sigma=\prod_{u\in\mu}x_{\sigma(u)}$.
Then the Macdonald polynomial
$\tilde{H}_\mu (x;q,t)$ have the following
explicit formula:
\begin{equation}
\tilde{H}_\mu (x;q,t)=
\sum_{\sigma:\mu\rightarrow\mathbb{Z}_{>0}}
q^{\mathrm{inv}(\sigma)}
t^{\mathrm{maj}(\sigma)} x^\sigma .
\end{equation}
\hfill$\square$
\end{theorem}

{}From the HHL formula,
we can show the following formula.
\begin{proposition}
For any partition $\mu$ and composition $\alpha$ of the same size,
one has
\begin{equation}\label{eq:HHL}
\sum_{p\in\mathcal{P}(\alpha)}
q^{\mathrm{inv}_\mu (p)}
t^{\mathrm{maj}_\mu (p)}=
\sum_{\eta\vdash |\mu|}
K_{\eta,\alpha}
\tilde{K}_{\eta,\mu}(q,t),
\end{equation}
where $\eta$ runs over all partitions of size $|\mu|$.
\end{proposition}
{\bf Proof.}
Indeed, \cite{HHL}, if two fillings $\sigma$ and $\sigma '$
belong to the same $s_\infty$ orbit, then
$\mathrm{Inv}(\sigma )=\mathrm{Inv} (\sigma ')$,
$\mathrm{Des}(\sigma )=\mathrm{Des} (\sigma ')$.
\hfill$\blacksquare$
\begin{corollary}
The (modified) Macdonald polynomial $\tilde{H}_\mu (x;q,t)$
have the following
expansion in terms of the monomial symmetric functions
$m_{\lambda}(x)$:
\begin{equation}
\tilde{H}_\mu (x;q,t)=
\sum_{\lambda\vdash |\mu|}
\left(
\sum_{p\in\mathcal{P}(\lambda )}
q^{\mathrm{inv}_\mu(p)}
t^{\mathrm{maj}_\mu(p)}
\right)
m_{\lambda}(x),
\end{equation}
where $\lambda$ runs over all partitions of size $|\mu|$.
\hfill$\square$
\end{corollary}

To find combinatorial interpretation of the
Kostka--Macdonald polynomials
$\tilde{K}_{\lambda,\mu}(q,t)$
remains significant open problem.
Among many important partial results
about this problem, we would like
to mention the following theorem
also due to Haglund--Haiman--Loehr:
\begin{theorem}[\cite{HHL} Proposition 9.2]
If $\mu_1\leq 2$, we have
\begin{equation}\label{eq:HHL_kosmac}
\tilde{K}_{\lambda,\mu}(q,t)=
\sum_{p\in\mathcal{P}_+(\lambda )}
q^{\mathrm{inv}_\mu (p)}
t^{\mathrm{maj}_\mu (p)}.
\end{equation}
\hfill$\square$
\end{theorem}
It is interesting to compare this formula with the
formula obtained by S. Fishel \cite{Fi},
see also \cite{Kir1},~\cite{KiSh}.

Concerning validity of the formula
Eq.(\ref{eq:HHL_kosmac}),
we state the following conjecture.
\begin{conjecture}\label{conj:HHL_kostka}
Explicit formula for
the Kostka--Macdonald polynomials
\begin{equation}
\tilde{K}_{\lambda,\mu}(q,t)=
\sum_{p\in\mathcal{P}_+(\lambda )}
q^{\mathrm{inv}_\mu (p)}
t^{\mathrm{maj}_\mu (p)}.
\end{equation}
is valid if and only if at least
one of the following two conditions
is satisfied.
\begin{enumerate}
\item[(i)]
$\mu_1\leq 3$ and $\mu_2\leq 2$.
\item[(ii)]
$\lambda$ is a hook shape.
\hfill$\square$
\end{enumerate}
\end{conjecture}
\begin{example}
As an example, the following is the list of
the tabloids associated with
the highest weight paths of weight $\lambda =(3,2,1)$
and shape $\mu =(4,2)$.
Here we also include the values of the
statistics in the form $(\mathrm{maj}_\mu,\mathrm{inv}_\mu):$
$$\begin{array}{|c|c|c|}\hline
\Yvcentermath1\young(2111,32)\,\,\,\,(2,3)&
\Yvcentermath1\young(2111,23)\,\,\,\,(1,4)&
\Yvcentermath1\young(1211,32)\,\,\,\,(1,3)
\begin{picture}(0,1.7)\end{picture}
\\[4mm]\hline
\Yvcentermath1\young(1211,23)\,\,\,\,(2,2)&
\Yvcentermath1\young(2211,31)\,\,\,\,(1,4)&
\Yvcentermath1\young(2211,13)\,\,\,\,(1,5)
\begin{picture}(0,1.7)\end{picture}
\\[4mm]\hline
\Yvcentermath1\young(3211,21)\,\,\,\,(0,6)&
\Yvcentermath1\young(3211,12)\,\,\,\,(0,5)&
\Yvcentermath1\young(1121,32)\,\,\,\,(2,2)
\begin{picture}(0,1.7)\end{picture}
\\[4mm]\hline
\Yvcentermath1\young(1121,23)\,\,\,\,(2,1)&
\Yvcentermath1\young(2121,31)\,\,\,\,(1,3)&
\Yvcentermath1\young(2121,13)\,\,\,\,(1,4)
\begin{picture}(0,1.7)\end{picture}
\\[4mm]\hline
\Yvcentermath1\young(3121,21)\,\,\,\,(0,5)&
\Yvcentermath1\young(3121,12)\,\,\,\,(1,4)&
\Yvcentermath1\young(1321,21)\,\,\,\,(1,3)
\begin{picture}(0,1.7)\end{picture}
\\[4mm]\hline
\Yvcentermath1\young(1321,12)\,\,\,\,(0,4)&
\begin{picture}(0,1.7)\end{picture}&
\begin{picture}(0,1.7)\end{picture}
\\[4mm]\hline
\end{array}$$
Then the generating function is
\begin{align*}
\sum_{p\in\mathcal{P}_+(\lambda)}
q^{\mathrm{inv}_{\mu} (p)}
t^{\mathrm{maj}_{\mu} (p)}=
q^6+q^5 t+2 q^5+4 q^4 t+q^4+q^3 t^2+3 q^3 t+2 q^2 t^2+q t^2
\end{align*}
which is different from
\begin{align*}
\tilde{K}_{\lambda ,\mu}(q,t)=
q^6+q^5 t+2 q^5+3 q^4 t+q^4+q^3 t^2+3 q^3 t+2 q^2 t^2+q^2 t+q t^2.
\end{align*}
Even if we consider the special value $t=1,$
these two polynomials are distinct.
Yet other 
examples which show that the formula (\ref{eq:HHL_kosmac})
does not hold if the condition 
(i) and (ii) of Conjecture \ref{conj:HHL_kostka} break down,
is given in Example \ref{ex:222_42}.
Let us remark that the choice $\lambda =(2,2,2)$ and
$\mu=(3,3)$ will give an example of both specializations
$q=1$ and $t=1$ give distinct polynomials.
\hfill$\square$
\end{example}

\subsection{Generating function of tau functions}
\label{sec:sumoftau}
Our main result is an elementary proof for special case
$t=1$ of the formula Eq.(\ref{eq:HHL})
in the following form.

\begin{theorem}\label{th:main}
Let $\alpha$ be a composition and $\mu$ be a partition of
the same size. Then,
\begin{equation}\label{eq:main1}
\sum_{p\in\mathcal{P}(\alpha)}
q^{\mathrm{maj}_{\mu} (p)}=
\sum_{\eta\vdash |\mu|}K_{\eta,\alpha}~K_{\eta,\mu}(q,1).
\end{equation}
\hfill$\square$
\end{theorem}

\begin{conjecture}\label{con:main}
Let $\alpha$ be a composition and $\mu$ be a partition of
the same size. Then,
\begin{equation}\label{eq:main}
q^{-\sum_{i\geq 2}\alpha_i}
\sum_{p\in\mathcal{P}(\alpha)}
q^{\tau_\mu (p)}=
\sum_{\eta\vdash |\mu|}K_{\eta,\alpha}\tilde{K}_{\eta,\mu}(q,1).
\end{equation}
\hfill$\square$
\end{conjecture}
Here $\sum_{i\geq 2}\alpha_i$
is equal to the number of letters other than 1
contained in each path $p\in\mathcal{P}(\alpha )$.

Let us remark that in view of general definition of $\tau^{r,s}$,
our $\tau_\mu$ is related with $\tau^{1,1}$,
whereas $\mathrm{maj}_\mu$ is related with $\tau^{r,1}$
where $r$ is bigger than the length of weight $\alpha$
(see Section 6).
As for intermediate $\tau^{r,s}$,
see Conjecture \ref{conj:tau} for some further information.

\begin{example}
Let us consider case $\alpha =(4,1,1)$
and $\mu =(4,2)$.
The following is a list of paths $p$ and the
corresponding value of tau function
$\tau_{(4,2)}(p)$.
For example, the top left corner
\fbox{$111123$\hspace{3mm}$3$} means
$p=\fbox{1}\otimes\fbox{1}\otimes\fbox{1}\otimes
\fbox{1}\otimes\fbox{2}\otimes\fbox{3}$
and $\tau_{(4,2)}(p)=3$.
$$\begin{array}{|cc|cc|cc|cc|cc|cc|}
\hline
111123& 3&
111132& 2&
111213& 2&
111231& 3&
111312& 2&
111321& 3\\\hline
112113& 3&
112131& 4&
112311& 3&
113112& 3&
113121& 4&
113211& 2\\\hline
121113& 4&
121131& 5&
121311& 4&
123111& 5&
131112& 4&
131121& 5\\\hline
131211& 4&
132111& 3&
211113& 5&
211131& 6&
211311& 5&
213111& 6\\\hline
231111& 7&
311112& 5&
311121& 6&
311211& 5&
312111& 6&
321111& 4\\\hline
\end{array}$$
Summing up, LHS of Eq.(\ref{eq:main}) is
$$q^{-2}\sum_{p\in\mathcal{P}((4,1,1))}
q^{\tau_{(4,2)} (p)}=
q^5+4 q^4+7 q^3+7 q^2+7 q+4.$$
In the RHS of Eq.(\ref{eq:main}),
nontrivial contributions come from the following
4 terms:
\begin{eqnarray*}
&&
K_{(6),(4,1,1)}\tilde{K}_{(6),(4,2)}(q,t)+
K_{(5,1),(4,1,1)}\tilde{K}_{(5,1),(4,2)}(q,t)\\
&&
+K_{(4,2),(4,1,1)}\tilde{K}_{(4,2),(4,2)}(q,t)+
K_{(4,1,1),(4,1,1)}\tilde{K}_{(4,1,1),(4,2)}(q,t)\\
&=&1\cdot (1)+2\cdot (q^3+q^2+qt +q+t)\\
&&+1\cdot (2 q^4+q^3t +q^3+2 q^2 t+q^2+qt +t^2)\\
&&+1\cdot (q^5+q^4t +q^4+2 q^3 t+q^3+2 q^2 t+qt^2 +qt )\\
   &=&q^5+(t+3) q^4+(3 t+4) q^3+(4 t+3) q^2+\left(t^2+4 t+2\right)
   q+t^2+2 t+1.
\end{eqnarray*}
By setting $t=1$, we get Eq.(\ref{eq:main}).
\hfill$\square$
\end{example}

\bigskip

\noindent
{\bf Proof of Theorem \ref{th:main}}

\begin{definition}
Let $\mu$ be a composition and $T$ be a tabloid of
size $|\mu |$.
Denote by $T^{(i)}$ the part of $T$ which is filled by numbers
from the interval $[\mu_{[i-1]+1},\mu_{[i]}]$.
Then
\begin{equation}
c_{\mu}(T)=\sum_{i\geq 1}c(T^{(i)}).
\end{equation}
\hfill$\square$
\end{definition}

\begin{lemma}
\begin{equation}
\sum_\eta K_{\eta ,\alpha}
K_{\eta ,\mu}(1,t)=
\sum_T~t^{c_\mu (T)},
\end{equation}
where the sum in the right hand side runs over standard tabloids $T$ of shape $\alpha$.
\end{lemma}

\noindent
{\bf Proof.}
Recall the following three formulas from \cite{M}, Chapter VI.

\noindent
{\it Formula 1.}
\begin{equation}
K_{\lambda,\mu}(q,t)=K_{\lambda,\mu'}(t^{-1},q^{-1})
t^{n(\mu)}q^{n(\mu')},
\end{equation}
and thus
\begin{equation}
\sum_{\eta}K_{\eta,\lambda}K_{\eta,\mu}(q,t)=
\sum_{\eta}K_{\eta,\lambda}K_{\eta,\mu'}(t^{-1},q^{-1})
t^{n(\mu)}q^{n(\mu')}.
\end{equation}

As a corollary of the formulas above,
\begin{equation}
\sum_{\eta}K_{\eta,\lambda}K_{\eta,\mu}(q,1)=
\sum_\eta K_{\eta,\lambda}K_{\eta,\mu'}(1,q^{-1})q^{n(\mu')}.
\end{equation}

\noindent
{\it Formula 2.}
\begin{equation}
J_\mu (x;1,t)=(t,t)_{\mu'} e_{\mu '}(x).
\end{equation}

\noindent
{\it Formula 3.}
\begin{equation}
(t,t)_r e_r(x)=
\sum_{\lambda\vdash r}
\frac{t^{n(\lambda ')}(t,t)_r}{H_\lambda (t)}S_\lambda (x,t)=
\sum_{\lambda\vdash r}K_{\lambda ,(1^r)}(t)S_\lambda (x,t).
\end{equation}

Therefore,
\begin{equation}
J_\mu (x;1,t)=\prod_{i\geq 1}
\left(
\sum_{\lambda^{(i)}\vdash\mu_i'}
K_{\lambda^{(i)},(1^{\mu_i'})}(t) S_{\lambda^{(i)}}(x,t)
\right),
\end{equation}
and after the plethystic change of variables
$X\longmapsto\frac{X}{1-t}$,
we obtain
\begin{equation}
\tilde{J}_\mu(x;1,t)=
\prod_{i\geq 1}\left(
\sum_{\lambda^{(i)}\vdash\mu_i}K_{\lambda^{(i)},(1^{\mu_i'})}(t)
s_{\lambda^{(i)}}(x)
\right)
\end{equation}
\hfill$\blacksquare$
\begin{claim}
\begin{equation}
\sum_{\lambda^{(i)}\vdash\mu_i'}K_{\lambda^{(i)},(1^{\mu_i'})}(t)
s_{\lambda^{(i)}}(x)=
\sum_{T}t^{\bar{c}(T)}~x^{Sh(T)},
\end{equation}
where the second sum runs over all standard tabloids $T$
of the size $r$, and $\tilde{c}(T)$ denotes either the charge of $T$,
or the value of tau function on the path corresponding to tabloid.
\end{claim}
In the case $\tilde{c}(T)=c(T)$ the charge of tabloid $T$,
this result is due to Lascoux--Sch\"{u}tzenberger;
in the case of $\tilde{c}(T)=n(\mu ')-\tau(T)$
this statement is a corollary of Theorem 7.4 and Corollary 6.13
from \cite{KSY}, where identification of tau function and
cocharge is given.
\begin{corollary}
\begin{equation}
\tilde{J}_\mu (x;1,t)=
\sum_{\lambda\vdash |\mu |}
\sum_T t^{\bar{c}_{\mu'} (T)}m_\lambda (x),
\end{equation}
where the second sum runs over the set of standard tabloids
of shape $\lambda$, and $\bar{c}_{\mu'}(T)=\sum \bar{c}(T^{(i)})$.
\end{corollary}

On the other hand,
\begin{eqnarray}
\tilde{J}_\mu (x;1,t)&=&
\sum_\eta K_{\eta,\mu}(1,t)s_\eta (x)
\nonumber\\
&=&\sum_\eta K_{\eta,\mu}(1,t)\sum_\lambda K_{\eta,\lambda}
m_{\lambda}(x)
\nonumber\\
&=&\sum_{\lambda}
\left(
\sum_\eta K_{\eta,\mu}(1,t)K_{\eta,\lambda}
\right) m_\lambda(x),
\end{eqnarray}
and therefore
\begin{equation}
\sum_{\eta}K_{\eta,\mu}(1,t)K_{\eta,\lambda}=
\sum_{T}t^{\tilde{c}_{\mu'}(T)},
\end{equation}
where the sum in the right hand side runs over the set of
standard $\lambda$-tabloids.

Finally,
\begin{equation}\label{eq:proof:a}
\sum_{\eta}~K_{\eta,\alpha}K_{\eta,\mu}(q,1)=
\sum_{\eta}~K_{\eta,\alpha}K_{\eta,\mu'}(1,q^{-1})~q^{n(\mu')}=
\sum_{T}~q^{n(\mu ')-\tilde{c}_{\mu}(T)},
\end{equation}
where the the third sum runs over the set of all standard $\alpha$-tabloids.

It remains to observe that according to Lemma~3.4,
$$\sum_{T}~{n(\mu ')-\tilde{c}_{\mu}(T)}= \sum_{i}~\biggl({\mu_{i} \choose 2}-
c(p_{[i]}) \biggr)=\tau_{\mu}-\sum_{i}~\mu_{i}~(1-\delta_{1,p_{1}^{[i]}})
=\mathrm{maj}_{\mu} (p).$$
\hfill$\blacksquare$

\subsection{Comments on generalizations of
Section \ref{sec:sumoftau}}
In order to clarify nature of tau statistics,
we consider possible generalizations of
the results in Section \ref{sec:sumoftau}.

\subsubsection{Regularization map and parabolic Kostka polynomials}
The main objective of this Section is
to give an interpretation of the energy 
statistics partition function for the BBS
as the value of a certain parabolic Kostka 
polynomial, see e.g., \cite{Kir,KiSh}.
This observation allows to write
a fermionic formula for the parabolic 
Kostka polynomials in question, see e.g. \cite{Kir},
as well as appears to be useful in the study of the BBS,
see e.g. \cite{KSY}.
\begin{definition}
Let $p$ be a path of type $\bigotimes_i B^{r_i,s_i}$ and
weight $\lambda$, define regularization $\tilde{p}=\mathrm{reg}(p)$
of the path $p$ to be
\begin{equation}
\mathrm{reg}(p)=(1\cdots n-1)^{\lambda_n}\cdots
(1\cdots i)^{\lambda_{i+1}}\cdots 1^{\lambda_2}p,
\end{equation}
where $(1\cdots i):=1\otimes\cdots\otimes i$,
and we have omitted all symbols $\otimes$.
\hfill$\square$
\end{definition}
Let $\tilde{T}$ be semi-standard Young tableau (i.e., highest weight)
corresponding to the regularized path $\tilde{p}$.

\begin{lemma}[\cite{KSY}, Lemma 7.2]
Assume that all $r_i=1$, then
\begin{equation}
\tau (p)=\tau (\tilde{p})+\mathrm{const}
=\bar{c}(\tilde{T})+\mathrm{Const,}
\end{equation}
where 
$$\mathrm{Const}=L(L-\mu_{1})-{L \choose 2}+\sum_{a=1}^{n}a{\mu_{a} \choose 2}+
\sum_{1 \le a < b \le n}a\mu_{a}\mu_{b},$$
and $L=\sum_{a}~m_{a}.$
\hfill$\square$
\end{lemma}
\begin{example}
Let $p_1=3332221$, then
$$\tilde{p_1}=1212121113332221.$$
We have $\tau (p_1)=7$, $\tau (\tilde{p}_1)=46$,
and $\tau (\tilde{p}_1)-\tau (p_1)=39$.
On the other hand, let  $p_2=3223123$, then
$$\tilde{p_2}=1212121113223123,$$
and $\tau (p_2)=14$, $\tau (\tilde{p}_2)=53$,
so that $\tau (\tilde{p}_2)-\tau (p_2)=39$,
as expected.
Time evolution of $p_2$ is
\begin{center}
\fbox{$
\begin{array}{lllllll}
 3 & 2 & 2 & 3 & 1 & 2 & 3 \\
 1 & 1 & 1 & 2 & 3 & 1 & 2 \\
 1 & 1 & 1 & 1 & 2 & 3 & 1 \\
 1 & 1 & 1 & 1 & 1 & 2 & 3 \\
 1 & 1 & 1 & 1 & 1 & 1 & 2 \\
 1 & 1 & 1 & 1 & 1 & 1 & 1 \\
 1 & 1 & 1 & 1 & 1 & 1 & 1
\end{array}
$}
\end{center}
\hfill$\square$
\end{example}
\begin{corollary}
\begin{equation}
\sum_\eta K_{\eta ,\mu}K_{\eta ,\lambda}(q)=
q^{C_1}K_{\Lambda ,(\kappa ,\lambda)}(q^{-1}),
\end{equation}
where $C_1$ is a constant and partitions $\Lambda$
and $\kappa$ are defined as follows;
$\Lambda =(\Lambda_1,\cdots,\Lambda_n)$, and
$\Lambda_i=\sum_{a\geq i}\mu_a$, $i=1,\cdots,n$;
$\kappa =(\kappa_1,\cdots,\kappa_n)$,
$\kappa_i=\sum_{a\geq i+1}$,
$i=1,\cdots,n-1$.
\end{corollary}
{\bf Proof.}
Indeed, we have
\begin{equation}
\mathrm{LHS}=\sum_{p\in\mathcal{P}(\lambda )}q^{\tau_\mu (p)}
=q^{C_2}\sum_{\tilde{T}}q^{\bar{c}(\tilde{T})}=
q^{C_3}K_{\Lambda,(\kappa ,\lambda)}(q^{-1}),
\end{equation}
where summation in the third term runs over all
Littlewood--Richardson tableaux of shape $\Lambda$ and weight
$(\kappa ,\lambda)$.
\hfill$\blacksquare$

\bigskip

\begin{conjecture}
Let $R_1=({s_i}^{r_i})_{1,2,\cdots}$ be a sequence
of rectangles, then
\begin{equation}\label{eq:conj:genmainth}
\sum_\eta K_{\eta ,R_1}K_{\eta ,R_2}(q)=
q^{C}K_{\Lambda ,(\kappa ,R_2)}(q^{-1}),
\end{equation}
where $\Lambda$ denotes partition $\left(\sum_{a\geq i}s_a\right)^{r_i}$,
$\kappa =\left(\sum_{a\geq i+1}s_a\right)^{r_i}$.
\hfill$\square$
\end{conjecture}

{\it It's} well-known that parabolic Kostka polynomials
$K_{\lambda ,R}(q)$, where $R=\{ (s_a^{r_a}\}_{a=1,2,\cdots}$
is a dominant (i.e., $s_1\geq s_2\geq\cdots$)
sequence of rectangular shape partitions, satisfy
the so-called duality theorem
\begin{equation}
K_{\lambda ,R}(q)=
K_{\lambda ',R'}(q^{-1})q^{n(R)},
\end{equation}
where $R'$ denotes a dominant rearrangement of the
sequence of rectangular shape partitions $\{(s_a^{r_a})\}_{a=1,2,\cdots}$
and
\begin{equation}
n(R)=\sum_{1\leq a<b}\min (r_a,r_b)
\min(s_a,s_b).
\end{equation}

\paragraph{Questions 1.}
Can one define a set $\mathcal{P}(R_1,R_2)$
of paths $\pi$ of rectangular type $R_2=\bigotimes_a (s_a^{r_a})$
and rectangular weight $R_1=\bigotimes_b(\mu_b^{\eta_b})$,
and energy function $e(\pi)$ such that
\begin{equation}
KK_{R_1,R_2}(q):=
\sum_\pi q^{e(\pi)}=
\sum_\eta K_{\eta ,R_1}K_{\eta ,R_2}(q).
\end{equation}
\hfill$\square$

\paragraph{Questions 2.}
Are there some ``physical interpretations"
of the duality theorems
\begin{equation}
K_{\lambda ,R}(q)=q^{n(R)}K_{\lambda ',R'}(q)
\end{equation}
and
\begin{equation}
KK_{R_1,R_2}(q)=q^{n(R_2)}KK_{R_1',R_2'}(q^{-1})
\end{equation}
for polynomials $K_{\lambda ,R}(q)$ and
$KK_{R_1,R_2}(q)$?
\hfill$\square$

\subsubsection{Generating functions of generalized tau statistics}
Let us consider generating functions of
the generalized tau statistics:
\begin{equation}
\tau^{r,s}(p):=\mathrm{maj}(u^{(r)}_s\otimes p),
\end{equation}
where $u^{(r)}_s \in B^{r,s}$ is the highest
element and $p\in B^{r_1,s_1}\otimes\cdots\otimes B^{r_L.s_L}$.
Let $\eta$ be a partition and $R$  be a sequence of rectangles, denote by
$K_{\eta, R}(q)$ the corresponding parabolic Kostka polynomial.
For any partition $\lambda$,
let $K_{\eta, \lambda}:=K_{\eta, \lambda}(1)$ be the 
corresponding Kostka number.
Then we have:
\begin{conjecture}\label{conj:tau}
Let $R=\{(r_i,s_i)_{i=1}^L\}$
be a sequence of rectangles.
Denote by $\mathcal{P}_R$ the set of all paths
(including non-highest elements)
corresponding to the tensor product of crystals
$B^{r_1,s_1}\otimes B^{r_2,s_2}\otimes\cdots\otimes
B^{r_L,s_L}$.
Then,
\begin{equation}\label{eq:conj:tau}
\sum_{p\in\mathcal{P}_R,\,
{\rm wt}(p)=\lambda}q^{\tau^{r,s}(p)}=
q^{\phi_r}\sum_\eta K_{\eta ,R}(q^{-1})K_{\eta ,\lambda},
\end{equation}
where $s\geq\max_i\{ s_i\}$, $\phi_r\in\mathbb{Z}$.
\hfill$\square$
\end{conjecture}
In other words, the statistics $\tau^{r,s}$ defines the essentially
unique class of polynomials although its definition
depends on choice of $u^{(r)}_s$.
As we see in the following example, this is
not obvious from definition of $\tau^{r,s}$.

\begin{example}
Let us consider the case $\lambda =(4, 3, 3, 2, 2)$
and $B^{3,2}\otimes B^{2,2}\otimes B^{2,2}$.
Then we have total of 759 paths, and by direct
computations, we have the following summation
over all 759 paths.
\begin{eqnarray*}
\sum q^{\tau^{1,s}(p)}&=&
q^{10}+8 q^9+33 q^8+89 q^7+161 q^6+198 q^5+163 q^4+82 q^3+24 q^2\\
\sum q^{\tau^{2,s}(p)}&=&
q^{13}+8 q^{12}+33 q^{11}+89 q^{10}+161 q^9+198 q^8+163 q^7+82
   q^6+24 q^5\\
\sum q^{\tau^{3,s}(p)}&=&
q^{12}+8 q^{11}+33 q^{10}+89 q^9+161 q^8+198 q^7+163 q^6+82 q^5+24
   q^4\\
\sum q^{\tau^{4,s}(p)}&=&
 q^{10}+8 q^9+33 q^8+89 q^7+161 q^6+198 q^5+163 q^4+82 q^3+24 q^2\\
\sum q^{\tau^{r,s}(p)}&=&
q^8+8 q^7+33 q^6+89 q^5+161 q^4+198 q^3+163 q^2+82 q+24,
\end{eqnarray*}
where $s=2,\cdots,5$ and,
in the last expression, $5\leq r\leq 10$.
However, if we look at specific paths,
for example,
$$b_1=\Yvcentermath1
\young(11,22,45)\otimes\young(23,34)\otimes\young(11,35)\,,
\qquad\qquad
b_2=\Yvcentermath1
\young(11,23,55)\otimes\young(12,34)\otimes\young(12,34)\,,$$
$$b_3=\Yvcentermath1
\young(12,24,35)\otimes\young(12,33)\otimes\young(11,45)\,,
\qquad\qquad
b_4=\Yvcentermath1
\young(22,33,44)\otimes\young(11,25)\otimes\young(11,35)\,.
$$
Then we have
$$\begin{array}{|c|cccccccccc|}
\hline
&\tau^{1,5}&\tau^{2,5}&\tau^{3,5}&\tau^{4,5}&
\tau^{5,5}&\tau^{6,5}&\tau^{7,5}&\tau^{8,5}&
\tau^{9,5}&\tau^{10,5}\\
\hline
b_1& 5 & 7 & 9 & 6 & 3 & 3 & 3 & 3 & 3 & 3 \\
b_2& 4 & 9 & 8 & 8 & 2 & 2 & 2 & 2 & 2 & 2 \\
b_3& 7 & 7 & 6 & 3 & 2 & 2 & 2 & 2 & 2 & 2 \\
b_4& 9 & 9 & 5 & 3 & 3 & 3 & 3 & 3 & 3 & 3\\\hline
\end{array}$$
In particular, dependences of $\tau^{r,s}(b)$ on $r$
are different for each $b$.
\hfill$\square$
\end{example}
\begin{remark}
We use the same notations of Conjecture \ref{conj:tau}.
Then, as we will see in Corollary \ref{cor:E=maj} below,
we have $\phi_r=\phi_{l(\lambda)}$ for all
$r\geq l(\lambda)$, where $l(\lambda)$
is length of weight $\lambda$.
\hfill$\square$
\end{remark}

\subsection{Generating functions related with the
energy statistics on the set of rectangular paths}
Let us consider generating function with respect to
$\bar{E}$ statistics.
This is a more traditional problem compared with
it for $\tau^{r,s}$.
Indeed, special cases of
this problem as well as its restriction to the set of
highest weight paths were considered by several
authors, see, e.g., \cite{NY,HKKOTY,O,Sch1, ScW}
and references therein.
\begin{conjecture}
Let $R=\{(r_i,s_i)_{i=1}^L\}$
is a sequence of rectangles.
Denote by $\mathcal{P}_R$ the set of all paths
corresponding to tensor product of crystals
$B^{r_1,s_1}\otimes B^{r_2,s_2}\otimes\cdots\otimes
B^{r_L,s_L}$.
Then,
\begin{equation}\label{eq:conj}
\sum_{p\in\mathcal{P}_R,\,
{\rm wt}(p)=\lambda}q^{\bar{E}(p)}=
\sum_{\eta\vdash |\lambda|}
K_{\eta ,R}(q)K_{\eta ,\lambda}.
\end{equation}
\hfill$\square$
\end{conjecture}
{\bf Comments.} An algebra--geometric definition of the parabolic Kostka 
polynomials (also known as  generalized Kostka polynomials) has been 
introduced in \cite{SW} (see also \cite{KiSh}) as a natural generalization of 
the well--known formula, \cite{M}, p.244,~(1), for the Kostka--Foulkes 
polynomials in terms of a $q$-analogue of the Kostant partition function. 
Based on the study of combinatorial properties of the algebraic Bethe ansatz, 
a fermionic formula for the parabolic Kostka polynomials has been discovered 
by the first author in the middle of 80's of the last century and has been 
proved in the full generality in \cite{KSS}. A ``path realization'' of the 
Kostka--Foulkes polynomials has been obtained in \cite{NY}, and finally, the 
formula
$$\sum_{p\in\mathcal{P}_{+,R},\,
{\rm wt}(p)}~q^{\bar{E}(p)}=K_{\lambda ,R}(q)$$
has been proved in \cite{ScW}.
\hfill$\square$
\begin{example}
Let us consider the case
$\lambda =(4,6,3,1)$ and
$R=\{(2,2),(2,2),(3,2)\}$, i.e.,
$B^{2,2}\otimes B^{2,2}\otimes B^{3,2}$.
Then we have the following nine paths:
$$
\begin{array}{|c|c||c|c|}
\hline
p&\bar{E}(p)&p&\bar{E}(p)\\\hline
\Yvcentermath1\young(11,22)\otimes\young(12,23)\otimes\young(12,23,34)&10&
\Yvcentermath1\young(11,22)\otimes\young(22,33)\otimes\young(11,22,34)
\begin{picture}(0,2)
\end{picture}&9
\\[15pt]\hline
\Yvcentermath1\young(11,22)\otimes\young(22,34)\otimes\young(11,22,33)&10&
\Yvcentermath1\young(12,23)\otimes\young(11,22)\otimes\young(12,23,34)
\begin{picture}(0,2)
\end{picture}&11
\\[15pt]\hline
\Yvcentermath1\young(12,23)\otimes\young(12,23)\otimes\young(11,22,34)&10&
\Yvcentermath1\young(12,23)\otimes\young(12,24)\otimes\young(11,22,33)
\begin{picture}(0,2)
\end{picture}&11
\\[15pt]\hline
\Yvcentermath1\young(12,24)\otimes\young(12,23)\otimes\young(11,22,33)&11&
\Yvcentermath1\young(22,33)\otimes\young(11,22)\otimes\young(11,22,34)
\begin{picture}(0,2)
\end{picture}&12
\\[15pt]\hline
\Yvcentermath1\young(22,34)\otimes\young(11,22)\otimes\young(11,22,33)
\begin{picture}(0,2)
\end{picture}&12&&
\\[15pt]\hline
\end{array}$$
Therefore, the LHS of Eq.(\ref{eq:conj}) is
$$q^9+3q^{10}+3q^{11}+2q^{12}.$$
On the other hand, non-zero contributions for
RHS of Eq.(\ref{eq:conj}) comes from
$$\begin{array}{|c||c|c|}
\hline
\eta & K_{\eta,R}(q) & K_{\eta,\lambda}(q)\\\hline\hline
(6,4,3,1) & q^9+q^{10}    & 1      \\\hline
(6,5,2,1) & q^{10}+q^{11} & q      \\\hline
(6,4,4)   & q^{10}        & q      \\\hline
(6,5,3)   & q^{11}        & q+q^2  \\\hline
(6,6,2)   & q^{12}        & q^2+q^3\\\hline
\end{array}$$
Summing up,
\begin{eqnarray*}
\mathrm{RHS}&=&
(q^9+q^{10})\cdot 1+(q^{10}+q^{11})\cdot 1+
(q^{10})\cdot 1+(q^{11})\cdot 2+q^{12}\cdot 2\\
&=&q^9+3q^{10}+3q^{11}+2q^{12},
\end{eqnarray*}
which coincides with the LHS.
\hfill$\square$
\end{example}

\section{Discussion: $\tau^{r,s}$ and $\bar{E}$}
So far in this paper, we have considered
several statistics including generalized
tau statistics $\tau^{r,s}$ and more
traditional one $\bar{E}$.
Let us investigate several further aspects of
these two statistics.
Our main results in this section are
(i) to show that $\bar{E}$ belong to the class
of statistics $\tau^{r,s}$ and
(ii) to show that $\tau^{r,s}$ stabilize
when we increase the value of $r$.
The following proposition will be a key property.
\begin{proposition}\label{prop:large_r}
Let $b_{i.j}$ be the integer at the
$i$-th row, $j$-th column of the tableau representation
of $b\in B^{r',s'}$, and let
the highest element of $B^{r,s}$
be $u^{r,s}$.
Then we have
\begin{equation}
H(u^{r,s}\otimes b)=0
\end{equation}
if $r\geq b_{r',s'}$ and $s\geq s'$.
\end{proposition}
Note that $b_{r',s'}$ is the largest integer
in the tableau representation of $b$.\bigskip

\noindent
{\bf Proof.}
According to the algorithm presented in
Proposition \ref{prop:shimozono},
we have to compute the insertion
$b\leftarrow row(u^{r,s})$.
This is worked out in Lemma \ref{lem:large_r} below,
and shape of the resulting tableau
coincides with the concatenation of two tableaux
$b$ and $u^{r,s}$.
Hence $H(u^{r,s}\otimes b)=0$
due to Proposition \ref{prop:shimozono}.
\hfill$\blacksquare$
\begin{lemma}\label{lem:large_r}
Under the same assumptions of Proposition
\ref{prop:large_r}, the insertion
\begin{equation}
b\leftarrow
\underbrace{rr\cdots r\mathstrut}_s\cdots
\underbrace{22\cdots 2\mathstrut}_s
\underbrace{11\cdots 1\mathstrut}_s.
\end{equation}
gives the concatenation of tableaux
$u^{r,s}$ and $b$, i.e.,
\begin{equation}\label{eq:concatenation}
\unitlength 1pt
\Yvcentermath1
\Yboxdim23pt
\newcommand{\cdotsr}{\begin{picture}(2,2)
\multiput(-4.5,4)(5.8,0){3}{\circle*{1}}
\end{picture}}
\newcommand{\bichiichi}{b_{1,1}}
\newcommand{\bichini}{b_{1,2}}
\newcommand{\bichis}{b_{1,s'}}
\newcommand{\bniichi}{b_{2,1}}
\newcommand{\bnini}{b_{2,2}}
\newcommand{\bnis}{b_{2,s'}}
\newcommand{\brichi}{b_{r',1}}
\newcommand{\brni}{b_{r',2}}
\newcommand{\brs}{b_{r',s'}}
\newcommand{\rp}{r'}
\newcommand{\bardelta}{\bar{\delta}}
\young(11\cdotsr 1\bichiichi\bichini\cdotsr\bichis,%
22\cdotsr 2\bniichi\bnini\cdotsr\bnis,%
\cdotsr\cdotsr\cdotsr\cdotsr\cdotsr\cdotsr\cdotsr\cdotsr,%
\rp\rp\cdotsr\rp\brichi\brni\cdotsr\brs,%
\cdotsr\cdotsr\cdotsr\cdotsr,%
rr\cdotsr r)\,\, .
\end{equation}
\end{lemma}
{\bf Proof.}
This insertion procedure can be divided
into two steps.
First, we show
\begin{align}\label{eq:proof0}
\tilde{b}:=(b\leftarrow
\underbrace{rr\cdots r\mathstrut}_s\cdots
\underbrace{\bar{\delta}\bar{\delta}\cdots
\bar{\delta}\mathstrut}_s
\underbrace{\delta\delta\cdots \delta\mathstrut}_s)
=&\,\,
\unitlength 1pt
\Yvcentermath1
\Yboxdim23pt
\newcommand{\cdotsr}{\begin{picture}(2,2)
\multiput(-4.5,4)(5.8,0){3}{\circle*{1}}
\end{picture}}
\newcommand{\bichiichi}{b_{1,1}}
\newcommand{\bichini}{b_{1,2}}
\newcommand{\bichis}{b_{1,s'}}
\newcommand{\bniichi}{b_{2,1}}
\newcommand{\bnini}{b_{2,2}}
\newcommand{\bnis}{b_{2,s'}}
\newcommand{\brichi}{b_{r',1}}
\newcommand{\brni}{b_{r',2}}
\newcommand{\brs}{b_{r',s'}}
\newcommand{\bardelta}{\bar{\delta}}
\young(\bichiichi\bichini\cdotsr\bichis%
\delta\delta\cdotsr\delta,%
\bniichi\bnini\cdotsr\bnis
\bardelta\bardelta\cdotsr\bardelta,%
\cdotsr\cdotsr\cdotsr\cdotsr\cdotsr\cdotsr\cdotsr\cdotsr,%
\brichi\brni\cdotsr\brs rrrr)\\[-12pt]
&
\underbrace{\hspace{92pt}}_{s'}
\underbrace{\hspace{92pt}}_{s}
\nonumber
\end{align}
where $\delta=r-r'+1$ and $\bar{\delta}=\delta +1$.
Note that we have $\delta\geq 1$
since (i) by the semi-standard property of $b$,
we have $b_{r',s'}\geq r'$ and (ii) by the assumption
$r\geq b_{r',s'}$, thus $r-r'\geq 0$.
Again from the assumption $r\geq b_{r',s'}$,
we have
\begin{equation}\label{eq:proof1}
r-i\geq b_{r'-i,s'}\qquad
(0\leq\forall i<r')
\end{equation}
by the semi-standard property $b_{k-1,s'}<b_{k,s'}$.
Consider the insertion $(b\leftarrow rr\cdots r)$.
{}From Eq.(\ref{eq:proof1}), we have
$r\geq r-(r'-1)\geq b_{r'-(r'-1),s'}=b_{1,s'}$.
Thus the first row of $(b\leftarrow rr\cdots r)$
is $b_{1,1}b_{1,2}\cdots b_{1,s'}rr\cdots r$,
and the remaining rows are identical to the corresponding
rows of $b$.
If $r=1$, this finishes the proof
(i.e., by $r\geq b_{r',s'}$ we have $b_{r',s'}=1$), therefore
let us consider the case $r>1$.
Assume that we have, for some $r\geq k>\delta$
and $\bar{k}=k+1$,
\begin{equation}
b^{\circ}:=
(b\leftarrow
\underbrace{rr\cdots r\mathstrut}_s\cdots
\underbrace{\bar{k}\bar{k}\cdots \bar{k}\mathstrut}_s
\underbrace{kk\cdots k\mathstrut}_s)
=
\unitlength 1pt
\Yvcentermath1
\Yboxdim23pt
\newcommand{\cdotsr}{\begin{picture}(2,2)
\multiput(-4.5,4)(5.8,0){3}{\circle*{1}}
\end{picture}}
\newcommand{\bichiichi}{b_{1,1}}
\newcommand{\bichini}{b_{1,2}}
\newcommand{\bichis}{b_{1,s'}}
\newcommand{\bniichi}{b_{2,1}}
\newcommand{\bnini}{b_{2,2}}
\newcommand{\bastichi}{b_{\ast ,1}}
\newcommand{\bastni}{b_{\ast ,2}}
\newcommand{\basts}{b_{\ast ,s}}
\newcommand{\bnis}{b_{2,s'}}
\newcommand{\brichi}{b_{r',1}}
\newcommand{\brni}{b_{r',2}}
\newcommand{\brs}{b_{r',s'}}
\newcommand{\bark}{\bar{k}}
\young(\bichiichi\bichini\cdotsr\bichis kk\cdotsr k,%
\bniichi\bnini\cdotsr\bnis
\bark\bark\cdotsr\bark,%
\cdotsr\cdotsr\cdotsr\cdotsr\cdotsr\cdotsr\cdotsr\cdotsr,%
\bastichi\bastni\cdotsr\basts rr\cdots r,%
\cdotsr\cdotsr\cdotsr\cdotsr,%
\brichi\brni\cdotsr\brs)\, ,
\end{equation}
where $\ast =r'-k+1$.
We insert $(k-1)$ into this $b^{\circ}$.
{}From Eq.(\ref{eq:proof1}) and the assumption
$k>\delta$, we have
$k>\delta =r-(r'-1)\geq b_{r'-(r'-1),s'}=b_{1,s'}$,
thus $k-1\geq b_{1,s'}$.
Therefore inserted $(k-1)$ bumps $k$ at $(s'+1)$-th
column of $b^{\circ}$.
As the next step, we have to insert $k$
to the second row of $b^{\circ}$.
By the similar reasoning, it bumps $(k+1)$
at $(s'+1)$-th column.
In this way, insertion of $(k-1)$ causes
downward shift of $(s'+1)$-th column of $b^{\circ}$
and addition of $(k-1)$ to the first row
of $(s'+1)$-th column.
Similarly, we see that the second insertion of $(k-1)$
causes shift and addition of $(k-1)$ to
the $(s'+2)$-th column.
Continuing in this way, we see that
successive insertions $(b^{\circ}\leftarrow (k-1)
(k-1)\cdots (k-1))$ gives the tableau
$b^{\circ}$ with replacement $k$ by $k-1$.
By induction, we can show Eq.(\ref{eq:proof0}).

Next, we consider the insertion of
$(\delta -1)^s=
(\delta -1)(\delta -1)\cdots (\delta -1)$
into $\tilde{b}$:
\begin{equation}
\tilde{\tilde{b}}:=(\tilde{b}\leftarrow
(\delta -1)^s).
\end{equation}
We denote the $i$-th row of $\tilde{b}$
(resp. $\tilde{\tilde{b}}$)
by $\tilde{b}_i$ (resp. $\tilde{\tilde{b}}_i$).
Although we independently repeat insertions of $(\delta -1)$
for $s$ times, we can argue more systematically
as follows.
Let us consider the first row of $\tilde{b}$:
\begin{equation}
\tilde{b}_1=b_{1,1}b_{1,2}\cdots b_{1,s'}
\underbrace{\delta\delta\cdots\delta}_{s}.
\end{equation}
As we saw in the last paragraph,
we have $b_{1,s'}\leq\delta$ from Eq.(\ref{eq:proof1}).
Suppose there are $k_1$ letters $\delta$
in the first row of $b$, i.e.,
$b_{1,s'-k_1+1}=b_{1,s'-k_1+2}=\cdots =b_{1,s'}=\delta$.
Thus there are $(s+k_1)$ letters $\delta$ in
the first row of $\tilde{b}$.
After inserting $(\delta -1)^s$,
precisely $s$ letters $\delta$ are bumped
and go down to the second row,
and the first row becomes
\begin{equation}
\tilde{\tilde{b}}_1=
b_{1,1}b_{1,2}\cdots b_{1,s'-k_1}
\underbrace{(\delta -1)(\delta -1)\cdots (\delta -1)}_{s}
\underbrace{\delta\delta\cdots\delta\mathstrut}_{k_1}.
\end{equation}
In particular, $k_1$ letters $\delta$ on the right of
the first row of $b$ are precisely
reproduced in the right part of $\tilde{\tilde{b}}_1$.

Now we have to consider the insertion of $s$ letters $\delta$
which are bumped from $\tilde{b}_1$
into the second row of $\tilde{b}$:
\begin{equation}
\tilde{b}_2=b_{2,1}b_{2,2}\cdots b_{2,s'}
\underbrace{(\delta +1)(\delta +1)\cdots (\delta +1)}_{s}.
\end{equation}
Again, from Eq.(\ref{eq:proof1}), we have $b_{2,s'}\leq\delta+1$,
and suppose that there are $k_2$ letters $(\delta +1)$
in the second row of $b$, i.e.,
$b_{2,s'-k_2+1}=b_{2,s'-k_2+2}=\cdots =b_{2,s'}=\delta+1$.
After inserting $\delta^s$,
precisely $s$ letters $(\delta+1)$ are
bumped and go down to the third row,
and the second row becomes
\begin{equation}
\tilde{\tilde{b}}_2=
b_{2,1}b_{2,2}\cdots b_{2,s'-k_2}
\underbrace{\delta\delta\cdots\delta\mathstrut}_{s}
\underbrace{(\delta +1)(\delta +1)\cdots (\delta +1)}_{k_2}.
\end{equation}
Again, $k_2$ letters $(\delta +1)$ on the right of
the second row of $b$ are precisely
reproduced in the right part of $\tilde{\tilde{b}}_2$.

As we see in the above discussions,
this procedure can be continued recursively.
Let the number of $\delta +i$ contained in
the $(i+1)$-th row of $b$ be $k_{i+1}$.
Then $\tilde{\tilde{b}}_{i+1}$ is
\begin{equation}
b_{i+1,1}b_{i+1,2}\cdots
b_{i+1,s'-k_{i+1}}
\underbrace{(\delta +i-1)(\delta +i-1)\cdots (\delta +i-1)}_s
\underbrace{(\delta +i)(\delta +i)\cdots (\delta +i)}_{k_{i+1}}
\end{equation}
and the right $k_{i+1}$ letters $(\delta +i)$ are
copy of those originally contained in
the right of $(i+1)$-th row of $b$.
In this way, $\tilde{\tilde{b}}$
contains copy of the letters in $b$.

As we have investigated insertion of $(\delta -1)^s$,
let us consider the insertion of
$(\delta -2)^s=
(\delta -2)(\delta -2)\cdots (\delta -2)$
into $\tilde{\tilde{b}}$.
Consider the row $\tilde{\tilde{b}}_1$.
Suppose that there are $l_1$ letters
$(\delta -1)$ within the first row of $b$, i.e.,
$b_{s'-(k_1+l_1)+1}=\cdots =b_{s'-k_1}=\delta -1$.
Thus there are total of $l_1+s$ letters $\delta -1$
in the row $\tilde{\tilde{b}}_1$.
Therefore, after insertion of $(\delta -2)^s$,
the row $\tilde{\tilde{b}}_1$ becomes
\begin{equation}
b_{1,1}b_{1,2}\cdots b_{1,s'-(k_1+l_1)}
\underbrace{(\delta -2)(\delta -2)\cdots (\delta -2)}_s
\underbrace{(\delta -1)(\delta -1)\cdots (\delta -1)}_{l_1}
\underbrace{\delta\delta\cdots\delta\mathstrut}_{k_1}.
\end{equation}
In particular, $l_1$ letters $(\delta -1)$
and $k_1$ letters $\delta$ are copy of the corresponding
letters of the first row of $b$.
As the result, $s$ letters $(\delta -1)$ are bumped
from $\tilde{\tilde{b}}_1$.
Therefore, for the row $\tilde{\tilde{b}}_2$,
we have to insert $(\delta -1)^s$,
and again obtain copy of the letters $\delta$
and $(\delta +1)$ contained in the second row of $b$.
We can continue this procedure until the bottom row
of $\tilde{\tilde{b}}$.
Therefore each row of the resulting tableau
$(\tilde{\tilde{b}}\leftarrow (\delta -2)^s)$
contains copy of at most two species of letters
in the original $b$.
We can recursively continue insertions of $(\delta -3)^s$,
$(\delta -4)^s$, $\cdots$, $1^s$,
and each insertion generates copy of part of letters of $b$.
Finally we get the result Eq.(\ref{eq:concatenation}).
\hfill$\blacksquare$

\bigskip

The following theorem shows that the
statistics $\bar{E}$ essentially
belong to the class of statistics $\tau^{r,s}$.
\begin{theorem}\label{th:E=tau}
Let $p=b_1\otimes b_2\otimes\cdots\otimes b_L
\in B^{r_1,s_1}\otimes B^{r_2,s_2}
\otimes\cdots\otimes B^{r_L,s_L}$,
and denote by $N$ the largest integer
contained in tableau representation of $p$.
If $r\geq N$ and $s\geq\max_i\{s_i\}$, we have
\begin{align}
\bar{E}(p)&=C-\tau^{r,s}(p),\\
C&=\sum_{i<j}\min{(r_i,r_j)}\cdot\min{(s_i,s_j)}.
\end{align}
\end{theorem}
{\bf Proof.}
Recall the definition of $\tau^{r,s}$
\begin{equation}
\tau^{r,s}(p)=\mathrm{maj}(u^{r,s}\otimes p),
\end{equation}
where maj is defined by Eq.(\ref{def:maj})
and $u^{r,s}$ is the highest element of $B^{r,s}$.
Within the definition of $\tau^{r,s}$,
let us first consider the energy functions
involving $u^{r,s}$ and next consider the
remaining ones.
As for the terms involving $u^{r,s}$, we have
\begin{equation}
H(u^{r,s}\otimes b_j^{(1)})=0\qquad
(1\leq\forall j\leq L)
\end{equation}
by Proposition \ref{prop:large_r}.
On the other hand, for the remaining contributions,
recall that
$b_i\in B^{r_i,s_i}$ and
$b_j^{(i+1)}\in B^{r_j,s_j}$.
Then from the definition of
normalizations, we have
\begin{equation}
\bar{H}(b_i\otimes b_j^{(i+1)})=
\min (r_i,r_j)\cdot\min (s_i,s_j)-
H(b_i\otimes b_j^{(i+1)}).
\end{equation}
Combining both contributions, we obtain the sought relation.
\hfill$\blacksquare$
\begin{corollary}\label{cor:E=maj}
Let $p=b_1\otimes b_2\otimes\cdots\otimes b_L
\in B^{r_1,s_1}\otimes B^{r_2,s_2}
\otimes\cdots\otimes B^{r_L,s_L}$,
and denote by $N$ the largest integer
contained in tableau representation of $p$,
and define $S=\max_i\{s_i\}$.
Then we have
\begin{align}
\tau^{N,S}(p)=\tau^{r,s}(p)
\end{align}
for all $r\geq N$ and $s\geq S$.
\hfill$\square$
\end{corollary}
Therefore Conjecture \ref{conj:tau}
means that we can define at most $N$
independent statistics $\tau^{r,s}(p)$,
where $N$ is the largest integer
contained in the path $p$,
and these statistics have essentially unique
generating function.

Let us remark physical interpretation of $\tau^{r,s}$.
For the paths (including non-highest elements)
of shape $B^{1,s_1}\otimes B^{1,s_2}\otimes\cdots
\otimes B^{1,s_L}$ and statistics $\tau^{1,S}$
($S=\max_i\{s_i\}$),
there is a straightforward
generalization of Theorem \ref{th:tau=rho},
see Section 4.1 of \cite{KSY}.
Under the same assumptions, $\tau^{1,S}$
is identified with cocharge of the
unrestricted rigged configurations.
It will be an interesting problem to find
a physical interpretation of the more general
$\tau^{r,s}(p)$.

\vspace{5mm}
\noindent
{\bf Acknowledgements:}
The work of RS is supported by
the Core Research for Evolutional Science and Technology
of Japan Science and Technology Agency.

\end{document}